\documentclass[a4paper,11pt]{paper}
\usepackage{amssymb,amsmath,amsthm}
\usepackage{array, hhline}
\usepackage{tikz}
\usepackage{authblk}

\usepackage[]{times}
\usepackage{fullpage}
\usepackage{graphicx}

\date{}

\usepackage{color}

\newcommand{\beginproof}{{\bf Proof}}
\newcommand{\finishedproof}{\hfill {\bf Q. E. D.}}

\newcommand{\ds}{\displaystyle}
\newcommand{\R}{\mathbb{R}}
\newcommand{\N}{\mathbb{N}}

\def\l2{L^2(\Omega)}
\def\ha1{H^1(\Omega)}
\def\h10{H^1_0(\Omega)}
\def\hdiv{H(\mathrm{div}; \Omega)}

\def\bun{\overline{{\bf u}}^n}

\def\opn{\overline{{p}}^n}
\def\othetan{\overline{\Theta}^n}
\def\oqn{\overline{{\bf q}}^n}
\def\oqk{\overline{{\bf q}}^k}

\def\bq{{\bf q}}
\def\bv{{\bf v}}

\def\bx{{\bf x}}
\def\bu{{\bf u}}
\def\bff{{\bf f}}

\newcommand{\ra}{\rangle}
\newcommand{\la}{\langle}

\newtheorem{proposition}{\sc \bf Proposition}[section]
\newtheorem{remark}{\bf Remark}[section]
\newtheorem{theorem}{\bf Theorem}[section]
\newtheorem{lemma}{\bf Lemma}[section]

\begin{document}

\title{A convergent mass conservative numerical scheme based on mixed finite elements for two-phase flow in porous media}

\author{ Florin A. Radu$^{1}$, Kundan Kumar$^{1}$, Jan M. Nordbotten$^{1}$, Iuliu S. Pop$^{1,2}$\\
{\small $^{1}$ Department of Mathematics, University of Bergen, P. O. Box 7800, N-5020 Bergen, Norway \\
{\small $^{2}$ Department of Mathematics and Computer Science, Eindhoven University of Technology,}\\
{\small P.~O. Box 513, 5600 MB Eindhoven, The Netherlands}\\
{\small e-mails: {\{florin.radu, jan.nordbotten, kundan.kumar\}@math.uib.no, i.pop@tue.nl}}}}

\maketitle




\begin{abstract}
{\bfseries Abstract.} In this work we present a mass conservative numerical scheme for two-phase flow in porous media. The model for flow consists on two fully coupled, non-linear equations: a degenerate parabolic equation and an elliptic equation. The proposed numerical scheme is based on backward Euler for the temporal discretization and mixed finite element method (MFEM) for the discretization in space. Continuous, semi-discrete (continuous in space) and fully discrete variational formulations are set up and the existence and uniqueness of solutions is discussed. Error estimates are presented to prove the convergence of the scheme. The non-linear systems within each time step are solved by a robust linearization method. This iterative method does not involve any regularization step. The convergence of the linearization scheme is rigorously proved under the assumption of a Lipschitz continuous saturation. Numerical results are presented to sustain the theoretical findings.
\end{abstract}


\section{Introduction}
\label{intro}
Two-phase porous media flow models are widely encountered in real-life applications of utmost societal relevance, including water and soil pollution, oil recovery, geological carbon dioxide sequestration, or nuclear waste management \cite{jan,rainer}. Such complex problems admit only in very simplified situations analytical solutions, therefore numerical methods for solving multiphase flow in porous media are playing a determining role in understanding and solving the problems.  Nevertheless, the design and analysis of robust, accurate and efficient numerical schemes is a very challenging task.

Here we discuss a numerical scheme for a two-phase porous media flow model. The fluids are assumed immiscible and incompressible and the solid matrix is non-deformable. The adopted formulation uses the global pressure and a complementary pressure, obtained by using the Kirchhoff transformation, as primary unknowns (see \cite{chavent,arbogast,chen}). This leads to a system of two coupled non-linear partial differential equations, a degenerate elliptic - parabolic one and an elliptic one.

Numerical methods for two-phase flow have been the object of intensive research in the last decades. The major challenge in developing efficient schemes is related to the degenerate nature of the problem. Due to this, the solution typically lacks regularity, which makes lower order finite elements or finite volumes a natural choice for the spatial discretization. In this respect, we refer to \cite{Fad,noc1} for Galerkin finite elements, to \cite{eymard2,ohlberger,michel} for finite volumes, to \cite{durlofsky,chen1997,chen2} for methods combining Galerkin finite elements combined with the mixed finite element method (MFEM), and to \cite{epsriv,riviere,karpinski} for the discontinuous Galerkin method. In all cases, the convergence of the numerical schemes is proved rigorously either by compactness arguments, or by obtaining {\it a priori} error estimates. {\it A posteriori} error estimates are obtained e.g. in \cite{CPV}. Furthermore, similar issues appear for the Richards equation, which is a simplified model for saturated/unsaturated flow in the case when the pressure of one phase is supposed to be constant. In this context we mention Galerkin finite elements \cite{noc1,soringeo}, MFEM based works \cite{arb_zhang, radu_SIAM, raduPhD, radu_NumMath,radu_wang,woodward,yotov}, multipointflux approximation (MPFA) \cite{klausen} or finite volume - MFEM combined methods \cite{eymard}. We also refer to \cite{radu_NMPDE,raduHabil,kumar1,kumar2} for related works concerning reactive transport in porous media. 

In this paper we propose a mass conservative scheme based on MFEM (lowest order Raviart-Thomas elements \cite{brezzi}) and backward Euler for numerical simulation of the two-phase flow in porous media. Continuous, semi-discrete (continuous in space) and fully discrete mixed variational formulations are defined. Existence and uniqueness of solutions is discussed, the equivalence with a conformal formulation being involved in the proof. We show the convergence of the numerical scheme and provide explicit order of convergence estimates. The analysis is inspired by similar results in \cite{arb_zhang,chen2,radu_SIAM,radu_NumMath}.

Typical problems involving flow in porous media, like e.g. water and soil pollution or nuclear waste management are spread over decades or even centuries, so that the use of relatively large time steps is a necessity. Due to this, implicit methods are a necessity (our choice here being the first-order backward Euler method, due to the low regularity of the considered problem).  Since the original model is non-linear, at each time step one needs to solve non-linear algebraic systems. In this work we propose a robust linearization scheme for the systems appearing at each time step, as a valuable alternative to modified Picard method \cite{celia} or Newton's method \cite{bergamaschi,park,radu2006,bastian,serge,radu2010,radu2011} or iterative IMPES \cite{kou1,kou2}.  Although the applicability of Newton's method for parabolic equations is well recognized, its convergence is not straightforward for degenerate equations, where the Jacobian might become singular. A possible way to overcome this is to regularize the problem. However, even in this case convergence is guaranteed only under a severe stability condition for the discretization parameters, see \cite{radu2006}. This has motivated the alternative, robust linearization scheme proposed in this work. The new scheme, called $L-$scheme from now on, does not involve the calculations of any derivatives and does not need a regularization step. The $L$-scheme combines the idea of a classical Picard method and the scheme presented in \cite{pop2004} for MFEM or \cite{yongpop,slod} for Galerkin finite elements. The $L$-scheme was proposed for two-phase flow in combination with the MPFA method in \cite{radu2015}, the proof of convergence there being only sketched and not made completely rigorous. 
We show here that the $L$-scheme for MFEM based discretizations converges linearly if the time step satisfies a mild condition. This robustness is the main advantage of the scheme when compared to the quadratic, but locally convergent Newton method.


Finally, we mention that the $L$-scheme can be interpreted as a non-linear preconditioner, because the linear systems to be solved within each iteration are much better conditioned than the corresponding systems in the case of modified Picard or Newton's method. We refer to \cite{list} for illustrative examples concerning the Richards equation, which is a particular case of the more general model considered in the present work.

To summarize, the main new contributions of this paper are
\begin{itemize}
\item{We present and analyze a MFEM based numerical scheme for two-phase flow in porous media. Order of convergence estimates are obtained.}
\item{We show the existence and uniqueness of the considered variational formulation. This is based on the equivalence between the conformal and the mixed formulations, which is proved here for the continuous and the time discrete models.}
\item{We present and analyze rigorously a robust, first-order convergent linearization method for MFEM based schemes for two-phase flow in porous media.}
\end{itemize}

The paper is structured as follows. In Section \ref{sec:model}, we present the model equations for two-phase flow in porous media and we define the discretization scheme. In Section \ref{sec:estimates} we analyze the convergence of the discretization scheme based on {\it a priori} error estimates. We also prove existence and uniqueness for the problem involved, and give stability estimates. A new MFEM linearization scheme is presented and analyzed in Section \ref{sec:lin}. Section \ref{sec:simulation} provides numerical examples confirming the theoretical results. The paper is ending with concluding remarks in Section \ref{sec:conclusions}.

\section{Mathematical model and discretization}\label{sec:model}
In this section we introduce the notations used in this work, the mathematical model and its MFEM/Euler implicit discretization (Problems $P$, $P^n$ and $P_h^n$). A linearization scheme (Problem $P_h^{n,i}$) is proposed to solve the  non-linear systems appearing at each time step.

The model is defined in the $d$-dimensional bounded domain $\Omega \subset {\mathbb R}^d$ having a Lipschitz continuous boundary $\Gamma$. Further $T > 0$ is the final computational time. We use common notations from functional analysis, e.g. $L^\infty(\Omega)$ is the space of essential bounded functions on $\Omega$, $L^2(\Omega)$ the space of square integrable functions on $\Omega$, or $H^1(\Omega)$  is the subspace of $L^2(\Omega)$ containing functions which have also first order derivatives in  $L^2(\Omega)$. We denote by $H^1_0(\Omega)$  the space of $H^1(\Omega)$ functions with a vanishing trace on $\Gamma$ and by $H^{-1}(\Omega)$ its dual. $\la \cdot,\cdot \ra$ denotes the inner product in $L^2(\Omega)$, or the duality pairing between $H^1_0(\Omega)$ and $H^{-1}(\Omega)$. Further, $\| \cdot \|$, $\| \cdot \|_1$ and $\| \cdot \|_{\infty}$ stand for the norms in $L^2(\Omega)$, $H^1(\Omega)$, respectively $L^\infty(\Omega)$. The functions in $\hdiv$ are vector valued having a $L^2$ divergence. The norm in $\hdiv$ is denoted by $\|\cdot \|_{div}$. $L^2(0, T; X)$ denotes the Bochner space of $X$-valued functions defined on $(0, T)$, where $X$ is a Banach space. Similarly, $C(0, T; X)$ are $X$-valued functions continuous (w. r. t. $X$ norm) on $[0, T]$. By $C$ we mean a generic positive constant, not depending on the unknowns or the discretization parameters and we denote by $L_f$ the Lipschitz constant of a (Lipschitz continuous) function $f(\cdot)$.\\

Further, we will denote by  $N \ge 1$ an integer giving the time step $\tau = {T/N}$. For a given $n \in \{1, 2, \dots, N\}$, the n$th$ time point is $t_n = n\tau$.  We will also use the following notation for the mean over a time interval. Given the function $g \in L^2(0, T; X)$ ($X$ being a Banach space like $L^2(\Omega)$, or $H^1(\Omega)$)), its time averaged over the interval $(t_{n-1}, t_n]$ is defined as
$$
\bar{g}^n := \frac 1 \tau \int_{t_{n-1}}^{t_n} g(t) dt.
$$
Clearly, this is an element in $X$ as well. 

The two-phase porous media flow model considered here assumes that the fluids are immiscible and incompressible, and that the solid matrix is non-deformable. By denoting with $\alpha = w, n$ the wetting and non-wetting phases, $s_\alpha, p_\alpha, \bq_\alpha, \rho_\alpha$ the saturation, pressure, flux and density of phase $\alpha$, respectively, the two-phase model under consideration reads (see e.g. \cite{bear,chavent,rainer,jan})
\begin{eqnarray}
\frac {\partial (\phi \rho_\alpha s_\alpha)}{\partial t} + \nabla \cdot  (\rho_\alpha \bq_\alpha) & = & 0, \quad \quad   \quad \quad   \, \, \, \, \alpha = w, n, \label{eq:2phase:1}\\[1ex]
 \bq_\alpha & = & - \dfrac{k_{r,\alpha}}{\mu_\alpha} k  (\nabla p_\alpha - \rho_\alpha {\bf g}), \quad \alpha = w, n,  \label{eq:2phase:2}\\[1ex]
s_w + s_n & = & 1, \label{eq:2phase:3}\\[1ex]
p_n - p_w & = &  p^{cap} (s_w),\label{eq:2phase:4}
\end{eqnarray}
where $ {\bf g}$ denotes the constant gravitational vector. Equation \eqref{eq:2phase:1} is a mass balance,  \eqref{eq:2phase:2} is the Darcy law, \eqref{eq:2phase:3} is an algebraic evidence expressing that all pores in the medium are filled by a mixture of the two fluid phases and \eqref{eq:2phase:4} is the capillary pressure relationship, with $p^{cap}(\cdot)$ supposed to be known. The porosity $\phi$, permeability $k$, the viscosities $\mu_\alpha$ are given constants and the relative permeabilities  $ k_{r,\alpha}(\cdot) $ are given functions of  $s_w$. We consider here a scalar permeability, but the results can be easily extended to the case when the permeability is positive-definite tensor.

In this paper we adopt a global/complementary pressure formulation \cite{arbogast,chavent,chen}. The global pressure (denoted by $p$) was introduced in \cite{chavent} and the complementary pressure in \cite{arbogast}. They are defined by
\begin{eqnarray}
p (\bx, s_w) := p_n (\bx)  - \int_0^{s_w} f_w(\bx, \xi) \dfrac{\partial p^{cap}}{\partial\xi}(\bx, \xi) d \, \xi, \\[2ex]
\Theta (\bx, s_w) := - \int_0^{s_w} f_w(\bx, \xi) \lambda_n(\bx, \xi)  \dfrac{\partial p^{cap}}{\partial\xi}(\bx, \xi) d \, \xi,\label{eq_theta}
\end{eqnarray}
where we denoted by $\lambda_\alpha := \dfrac{k_{r,\alpha}}{\mu_\alpha}$, $\alpha = w, n$ the phase mobilities and by $f_w := \dfrac{\lambda_w}{\lambda_w + \lambda_n}$ the fractional flow function. We note the use of Kirchhoff transformation above. In the new unknowns, the resulting system consists of two coupled non-linear partial differential equations, a degenerate parabolic one and an elliptic one. For more details on the modelling we refer to \cite{chen}, where the existence and uniqueness of a weak solution is proved for a Galerkin-MFEM formulation. In the new unknowns the system \eqref{eq:2phase:1}-\eqref{eq:2phase:4} becomes
\begin{eqnarray}
\partial_t s(\Theta) + \nabla \cdot \bq &=& 0, \label{eq_classic_1}\\[1ex]
\bq &=& - \nabla \Theta + f_w(s) \bu + \bff_1(s),  \label{eq_classic_2}\\[1ex]
\nabla \cdot \bu &=& f_2(s), \label{eq_classic_3}\\[1ex]
a(s) \bu &=& - \nabla p - \bff_3(s).  \label{eq_classic_4}
\end{eqnarray}
with $s := s_w$, $a(s) := \dfrac{1}{k \lambda(s)}$, $\bq$ the (wetting) flux, and $\bu$ the total flux. The equations hold true in $\Omega \times (0, T]$. The coefficient functions $s(\cdot), a(\cdot), f_w (\cdot), \bff_1(\cdot), f_2(\cdot), \bff_3(\cdot)$ are given and satisfy the assumptions listed below. The system is completed by initial conditions specified below, and by boundary conditions. For simplicity, we restrict our attention to homogeneous Dirichlet boundary conditions, but other kinds of conditions can be considered.

Also note that the results can be extended straightforwardly to the case when a source term $f_s$ is present on the right hand side of \eqref{eq_classic_1}. For the ease of presentation and in view of the analogy with the model considered in \cite{radu_NumMath}, the source terms are left out here. This simplifies the presentation of the convergence proof in Section \ref{sec:estimates}.

{\bf Problem $P$: Continuous mixed variational formulation.} \\
Find $\Theta, p  \in L^2(0, T; L^2(\Omega))$, $\bq \in L^2(0, T;(L^2(\Omega))^d)$,  $\bu \in L^2(0, T;\hdiv)$ such that there holds $s(\Theta) \in L^\infty(\Omega \times (0, T))$, $\int_0^t \bq(y)\, dy \in C(0, T;\hdiv) $, and 
\begin{eqnarray}
\la s(\Theta(t))-s(\Theta^0), w \ra + \la \nabla \cdot \int_0^t \bq(y) \, dy, w\ra &=& 0, \label{eq_cont_1}\\[1ex]
\la \int_0^t \bq(y) \, dy, \bv\ra - \la \int_0^t \Theta(y) \, dy,\nabla \cdot \bv \ra && \nonumber \\[1ex]
- \la \int_0^t  f_w(s(\Theta(y))) \bu (y)\,dy, \bv \ra
& = & \la \int_0^t \bff_1 (s(\Theta(y))) \,dy, \bv \ra, \label{eq_cont_2}  \\[1ex]
\la \nabla \cdot \bu(t), w \ra &=& \la f_2(s(\Theta(t))), w \ra, \label{eq_cont_3}\\[1ex]
\la a(s(\Theta(t))) \bu(t), \bv \ra - \la p(t), \nabla \cdot \bv \ra +  \la \bff_3(s(\Theta(t))), \bv \ra &=& 0 \label{eq_cont_4}
\end{eqnarray}
for all $t \in (0, T]$,  $w \in L^2(\Omega)$ and $\bv \in \hdiv$, with $\Theta(0) = \Theta_I \in L^2(\Omega)$.

The function $\Theta_I$ is given. For example, if the function $s(\cdot)$ is one to one, a natural initial condition is $\Theta_I = s^{-1}(s_I)$, where $s_I$ is a given initial saturation. By \eqref{eq_cont_1}, since $\int_0^t \bq(y) dy $ is continuous in time, it follows that $s(\Theta) \in C(0, T; L^2(\Omega))$. Then, since $s(\cdot)$ and $f_2(\cdot)$ are assumed continuous (see (A1) and (A3) below), from \eqref{eq_cont_3} one obtains that $\bu \in C(0. T; \hdiv)$. Similarly, $p$ is continuous in time as well, so \eqref{eq_cont_2}-\eqref{eq_cont_4} hold for all $t \in (0, T]$.

We now proceed with the time discretization for Problem $P$, which is achieved by the Euler implicit scheme. For a given $n \in \{1, 2, \dots, N\}$, we define the time discrete mixed variational problem at time $t_n$ (the time step is denoted by $\tau$):

{\bf Problem $P^n$: Semi-discrete variational formulation.} Let  $\Theta^{n-1}$ be given. Find $\Theta^n, p^n \in L^2(\Omega)$ and  $\bu^n,  \bq^n  \in \hdiv$ such that
\begin{eqnarray}
 \la s^n-s^{n-1},w \ra + \tau \la \nabla \cdot \bq^n, w \ra &=& 0, \label{eq_semi-discrete_1}\\ [1ex]
\la \bq^n,\bv \ra - \la \Theta^n,\nabla \cdot \bv \ra - \la f_w(s^n) \bu^n,\bv \ra &=& \la \bff_1(s^n), \bv \ra,\label{eq_semi-discrete_2}\\[1ex]
\la \nabla \cdot \bu^n, w \ra &=& \la f_2(s^n), w\ra, \label{eq_semi-discrete_3}\\[1ex]
 \la a(s^n) \bu^n, \bv \ra   - \la p^n, \nabla \cdot \bv \ra + \la \bff_3(s^n), \bv \ra &=& 0 \label{eq_semi-discrete_4}
\end{eqnarray}
for all $w \in L^2(\Omega)$, and $\bv \in \hdiv$. Initially we take $\Theta^0 = \Theta_I \in L^2(\Omega)$. Throughout this paper $s^k$ stands for $s(\Theta^k)$, $k \in \N$, making the presentation easier.

We can now proceed with the spatial discretization. For this let ${\mathcal{T}_h}$ be a regular decomposition of $\Omega \subset {\mathbb R}^d$ into closed $d$-simplices; $h$ stands for the mesh-size (see \cite{cia}). Here we assume $\overline\Omega = \cup_{T \in \mathcal T_h} T$, hence $\Omega$ is polygonal. Thus we neglect the errors caused by an approximation of a non-polygonal domain and avoid an excess of technicalities (a complete analysis in this sense can be found in \cite{noc1}).

The discrete subspaces $W_h \times V_h \subset L^2(\Omega) \times \hdiv$ are
defined as
\begin{equation}\label{discrete_spaces}
\begin{array}{l}
W_h := \{ p \in L^2(\Omega) \vert \ p \mbox{ is constant on each element }T \in \mathcal{T}_h \}, \\[2ex]
 V_h := \{ \bq \in \hdiv | \bq_{\vert T}({\bf x})= {\bf{a_T}} + b_T {\bf{x}}, {\bf{a_T}} \in \R^2, b_T \in \R \mbox{ for all } T \in \mathcal{T}_h \}.
\end{array}
\end{equation}
So $W_h$ denotes the space of piecewise constant functions, while $V_h$ is the $RT_0$ space (see \cite{brezzi}).

We will use the following projectors (see \cite{brezzi} and \cite{quarteroni}, p. 237):
\begin{equation}
P_h: L^2(\Omega) \rightarrow W_h,
\qquad \la P_h w - w, w_h \ra = 0, \label{eq_l2_projector}
\end{equation}
and 
\begin{equation}\label{eq_hdiv_projector}
\Pi_h: \hdiv \rightarrow V_h, \qquad
\la \nabla \cdot
(\Pi_h \bv - \bv) , w_h \ra = 0,
\end{equation}
for all $w \in \l2$,$\bv \in \hdiv$ and $w_h \in W_h$. For these operators we have
\begin{equation}\label{lemmaproiectors}
\| w - P_h w \| \leq C h \| w \|_{1}, \qquad \mbox{ respectively } \qquad 
\| \bv - \Pi_h \bv \| \leq C h \| \bv \|_{1}
\end{equation}
for any  $w \in H^1(\Omega)$ and $\bv \in (H^1(\Omega))^d$.

The fully discrete (non-linear) scheme can now be given. To simplify notations we use in the following the notation $s^{n}_h := s(\Theta^{n}_h)$, $n \in \N$.

{\bf Problem $P_h^{n}$: Fully discrete (non-linear) variational formulation.} Let $n \in \N, n \ge 1$,  and assume $\Theta^{n-1}_h$ is known. Find $\Theta^{n}_h,  p^{n}_h \in W_h$ and $\bq^{n}_h, \bu^{n}_h \in V_h$ such that there holds
\begin{eqnarray}
\la s^{n}_h - s^{n-1}_h, w_h \ra + \tau \la \nabla \cdot \bq^{n}_h, w_h \ra &=& 0, \label{eq1_discrete_nonlin}\\[1ex]
\la \bq^{n}_h, \bv_h \ra -\la \Theta^{n}_h,  \nabla  \cdot \bv_h \ra - \la f_w (s^{n}_h) \bu^{n}_h,  \bv_h \ra  &=& \la \bff_1 (s^{n}_h), \bv_h \ra, \label{eq2_discrete_nonlin}\\[1ex]
\la \nabla \cdot \bu^{n}_h, w_h \ra   &=& \la f_2(s^{n}_h), w_h \ra, \label{eq3_discrete_nonlin}\\[1ex]
\la a(s^{n}_h) \bu^{n}_h, \bv_h \ra - \la p^{n}_h, \nabla \cdot \bv_h \ra + \la \bff_3(s^{n}_h), \bv_h \ra&=& 0\label{eq4_discrete_nonlin}
\end{eqnarray}
for all $w_h \in W_h$ and all $\bv_h \in V_h$.

The fully discrete scheme (\ref{eq1_discrete_nonlin}) --  (\ref{eq4_discrete_nonlin}) is non-linear, and iterative schemes are required for solving it. Moreover, in (A1) it is only required that $s(\cdot)$ is monotone and Lipschitz, so it may have a vanishing derivative, which makes \eqref{eq_classic_1} degenerate. This is, indeed, the situation encountered in two-phase porous media flow models. In such cases, usual schemes such as the Newton method may not converge without performing a regularization step. This may affect the mass balance. For the Richards equation, this is proved in \cite{radu2006}. Following the ideas in \cite{pop2004,slod,yongpop}, we propose a robust, first order convergent linearization scheme for solving (\ref{eq1_discrete_nonlin}) -- (\ref{eq4_discrete_nonlin}).  The scheme is not requiring any regularization. Similar ideas can be applied in connection with any other spatial discretization method, see e.g. \cite{radu2015} for multipoint flux approximation. The analysis of the scheme is presented in Section \ref{sec:lin}.

Let $n \in N, n \ge 1$ be fixed. Assuming that $s(\cdot)$ is Lipschitz continuous, the following iterative scheme can be used to solve the non-linear problem (\ref{eq1_discrete_nonlin}) --  (\ref{eq4_discrete_nonlin}):

{\bf Problem $P_h^{n,i}$: Linearization scheme ($L$-scheme).} Let $L \ge L_s$, $i \in \N, i \ge 1$ and let  $\Theta^{n, i-1}_h \in W_h$ be given. Find $\Theta^{n,i}_h,  p^{n,i}_h \in W_h$ and $\bq^{n,i}_h, \bu^{n,i}_h \in V_h$ such that
\begin{eqnarray}
\la L (\Theta^{n,i}_h - \Theta^{n,i-1}_h) + s^{n,i-1}_h, w_h \ra  + \tau \la \nabla \cdot \bq^{n,i}_h, w_h \ra &=& \la s^{n-1}_h, w_h \ra,\label{eq1_discrete_lin}\\
\la \bq^{n,i}_h, \bv_h \ra -\la \Theta^{n,i}_h,  \nabla  \cdot \bv_h \ra - \la f_w (s^{n,i-1}_h) \bu^{n,i}_h,  \bv_h \ra &=&  \la \bff_1(s^{n,i-1}_h) , \bv_h \ra ,\label{eq2_discrete_lin}\\
\la \nabla \cdot \bu^{n,i}_h, w_h \ra   &=& \la f_2(s^{n,i-1}_h), w_h \ra, \label{eq3_discrete_lin}\\
\la a(s^{n,i-1}_h) \bu^{n,i}_h, \bv_h \ra - \la p^{n, i}_h, \nabla \cdot \bv_h \ra + \la \bff_3 (s^{n,i-1}_h), \bv_h \ra &=& 0\label{eq4_discrete_lin}
\end{eqnarray}
for all $w_h \in W_h$ and all $\bv_h \in V_h$. We use the notation $s^{n,i}_h := s(\Theta^{n,i}_h)$, $n \in \N$ and, as previously,  $s^{n}_h := s(\Theta^{n}_h)$, $n \in \N$. For starting the iterations, a natural choice is $\Theta^{n,0}_h := \Theta^{n-1}_h$ and, correspondingly, $s^{n,0}_h := s^{n-1}_h$. Note that since we prove below that the iterative scheme is a contraction, this choice is not compulsory for the convergence.

Throughout this paper we make the following assumptions:
\begin{itemize}
\item[(A1)]{The function $s(\cdot): \R \rightarrow [0, 1]$ is monotone increasing and Lipschitz continuous. }
\item[(A2)]{$a(\cdot)$ is Lipschitz continuous and there exists $a_\star, a^\star > 0 $ such that for all $y \in \R$ one has
\begin{equation}\label{bounds_a}
0 < a_\star \le a(y) \le a^\star < \infty.
\end{equation}
}
\item[(A3)]{$\bff_1(\cdot), f_2(\cdot), \bff_3(\cdot)$ and $f_w(\cdot)$ are Lipschitz continuous. Additionally,  $f_w(\cdot)$ is uniformly bounded.}
\item[(A4)]{There exits a constant $M_{\bu} < \infty$ such that $\| \bu \|_\infty \le M_{\bu}$, $\| \bu_h^n \|_\infty \le M_{\bu}$, and $\| \bu_h^{n, i} \|_\infty \le M_{\bu}$ for all $n \in \N$, the last two being uniformly in $h$ and $i$. Here $\bu$, $\bu_h^n$  and $\bu_h^{n, i}$ are the solution components in Problems P, P$_h^n$ and P$_h^{n,i}$ respectively. }
\item[(A5)]{The function $\Theta_I$ is in $L^2(\Omega)$.}
\end{itemize}

\begin{remark} \label{rem:A4}
The assumptions are satisfied in most situations of practical interest. Concerning (A4), for $\bu$ this is practically the outcome of the assumptions (A1) and (A3), which guarantee that for every $t \in [0, T]$ one has $f_2(s(\Theta(t))) \in L^\infty(\Omega)$, and that the $L^\infty$  norm is bounded uniformly w.r.t. time. Now, without being rigorous, we observe that by \eqref{eq_cont_3} one obtains $\bu(t) = - \nabla w(t)$, where $w$ satisfies  $- \Delta w(t) = f_2(s(\Theta(t)))$. Classical regularity theory (see e.g. \cite{Ladyz}, Thm. 15.1 in Chapter 3) guarantees that $\nabla w$ is continuous on the compact $\bar \Omega$, and that the $L^\infty$ norm can be bounded uniformly in time. For the approximation $\bu_h^n$, one can reason in the same manner, and observe that $\bu_h^n$ becomes the projection $\Pi_h (-\nabla w(t_n))$. Since $\|\nabla w(t_n)\|_\infty$ is bounded uniformly in time, the construction of the projector $\Pi_h$ (see e.~g. \cite{quarteroni}, Chapter 7.2) guarantees that $\bu_h^n$ satisfies the same bounds as $\nabla w$. Finally, case of $\bu_h^{n,i}$ is similar.  We also refer to \cite{radu_NumMath} for a similar situation but in the case of a one phase flow model, where conditions ensuring the validity of (A4) are provided.
\end{remark}

\begin{remark}\label{remark:extension_lipschitz}
One may relax the Lipschitz continuity for the parameter functions to a more general context, where e.~g. $f_w$ satisfies a growth condition
$|f_w(s_1) - f_w(s_2)|^2 \le C | \la f_w(s_1) - f_w(s_2), s_1- s_2 \ra | $. We do not exploit this possibility here, but refer to \cite{chen,radu_NumMath,radu_SIAM} for the details on the procedure.
\end{remark}

The following two technical lemmas will be used in Sections \ref{sec:estimates} and \ref{sec:lin}. Their proofs can be found e.g. in \cite{radu_NumMath} and \cite{thomas}, respectively.
\begin{lemma}\label{lemma_thomas_cont}  Given a $w \in L^2(\Omega)$, there exists a $\bv \in \hdiv$ such that
\begin{displaymath}
\nabla \cdot \bv = w  \,\,{ \rm{ and } }\,\,\|{\bv}\| \le C_\Omega \| w \|,
\end{displaymath}
with $C_\Omega > 0$ not depending on $w$.
\end{lemma}
\begin{lemma} \label{lemma_thomas} Given a $ w_h \in W_h$,  there exits a $ \bv_h \in V_h$  satisfying
\begin{displaymath}
\nabla \cdot \bv_h =  w_h  \,\,{ \rm{ and } }\,\,  \| \bv_h \| \le C_{\Omega,d}  \| w_h \|,
\end{displaymath}
with $C_{\Omega,d} > 0$ not depending on $w_h$ or mesh size.
\end{lemma}

Also, the following elementary result will be used
\begin{proposition}\label{prop_1}
Let ${\bf a}_k \in {\mathbb R}^d$ $(k \in \{1, \ldots, N\}, d \geq 1)$ be a set of $N$ vectors. It holds
\begin{eqnarray}
\sum_{n=1}^N \la {\bf a}_n, \sum_{k=1}^n {\bf a}_k \ra &=& \frac 1 2 \left\| \sum_{n=1}^N {\bf a}_n \right\|^2 + \frac 1 2 \sum_{n=1}^N \left\| {\bf a}_n \right\|^2\,.\label{lemma1_1}
\end{eqnarray}
\end{proposition}

\section{Analysis of the discretization: existence and uniqueness and { a priori} error estimates}\label{sec:estimates}

In this section we analyze the problems $P$, $P_n$ and $P^n_h$ introduced in Section \ref{sec:model}. The existence and uniqueness of a solution will be discussed in Subsection \ref{sec:equivalence}. For the continuous and semi-discrete cases this will be done by showing an equivalence with conformal variational formulations.   The convergence of the numerical scheme will be shown by deriving {\it a priori} error estimates. The main convergence result is given in Theorem \ref{th_main_result_estimates}. The convergence is established by assuming that the non-linear systems \eqref{eq1_discrete_nonlin} - \eqref{eq4_discrete_nonlin} are solved exactly. We refer to Section \ref{sec:lin} for the analysis of the linearization scheme \eqref{eq1_discrete_lin}-\eqref{eq4_discrete_lin}, which was proposed in the previous section to solve these non-linear algebraic systems numerically.

\subsection{Existence and uniqueness for the variational problems} \label{sec:equivalence}
In this subsection we discuss the existence and uniqueness of the continuous, semi-discrete, and fully discrete variational formulations for the considered model \eqref{eq_classic_1}-\eqref{eq_classic_4}. We establish an equivalence between the continuous mixed formulation and a conformal formulation, which will deliver the existence and uniqueness for the continuous case. The semi-discrete case can be treated analogously. For the fully discrete case we prove below the uniqueness. Existence can be proved by fixed point arguments, using e.g. Lemma 1.4, p. 164 in \cite{Temam}. We omit the details here as the existence is also a direct consequence of the results in Section \ref{sec:lin}, where the convergence of a linear iterative scheme is proved. The limit of this iteration is exactly a solution for the fully discrete system.

A conformal variational formulation for the model  \eqref{eq_classic_1}-\eqref{eq_classic_4} reads:

{\bf Problem $PC$: Continuous conformal variational formulation.} Let  $\Theta_I \in L^2(\Omega)$ be given. Find $\Theta_C$, $p_C  \in L^2(0, T; H^1_0(\Omega))$, such that $\partial_t s(\Theta_C) \in L^2(0, T; H^{-1}(\Omega))$, $\Theta_C(0) = \Theta_I $, and for all $v \in L^2(0, T; H^1_0(\Omega))$ and $w \in H^1_0(\Omega)$ one has
\begin{eqnarray}
\int_0^T \la \partial_t s(\Theta_C), v \ra dt + \int_0^T \big\la \nabla \Theta_C + \frac {f_w(s\Theta_C)}{a(s(\Theta_C))}\nabla p_C, \nabla v\big\ra dt && \nonumber \\
+ \int_0^T \big\la \frac {f_w(s\Theta_C)}{a(s(\Theta_C))}\bff_3(s(\Theta_C)) - \bff_1(s(\Theta_C)), \nabla v\big\ra dt &=& 0, \label{eq_cf_cont_1}\\[1ex]
\big \la \frac 1 {a(s(\Theta_C))} \big(\nabla p_C + \bff_3(s(\Theta_C))\big), \nabla w \big\ra &=& -\la f_2(s), w \ra . \label{eq_cf_cont_3}
\end{eqnarray} 

The existence and uniqueness of a solution for Problem $PC$ 
has been studied intensively in the past. Closest to the framework considered here is \cite{chen} (see also \cite{chen2}). There the existence and uniqueness is proved, however, for the case that the inverse of $s(\cdot)$ is Lipschitz (the so-called slow diffusion case). Here we assume $s(\cdot)$ Lipschitz but not necessarily strictly increasing, which is a fast diffusion case. Other relevant references for the existence and uniqueness are \cite{Alt_diB, Amaz, Fad, Kroener_Luckhaus}. Also, we refer to \cite{Cances_P} for the existence of a solution in heterogeneous media, where the phase pressure differences may become discontinuous at the interface separating two homogeneous blocks. 

Having this in mind, one can use the existence and uniqueness results for the conformal formulation to obtain the existence of a solution for Problem $P$ (as by-product also establish its regularity). The equivalence is established in Proposition \ref{proposition:equivalence}, whose proof follows the ideas in \cite{radu_SIAM}, Proposition 2.2.
\begin{proposition}\label{proposition:equivalence}
Let $\Theta_C, p_C  \in L^2(0, T; H^1_0(\Omega))$ be a solution to Problem $PC$, define $s_C = s(\Theta_C)$ and assume that (A1)-(A5) hold true. Then, a solution to Problem $P$ is given by $\Theta = \Theta_C$, $p = p_C$, $\bq = -\nabla \Theta_C -  \frac {f_w (s_C)} {a(s_C)}\left( \nabla p_C + \bff_3(s_C) \right) + \bff_1(s_C)$ and $\bu = - \frac 1 {a(s_C)}\left( \nabla p_C + \bff_3(s_C) \right)$. Conversely, if $(\Theta,\bq)  \in L^2(0, T; L^2(\Omega)) \times L^2(0, T;(L^2(\Omega))^d)$,  $(p, \bu) \in L^2(0, T;L^2(\Omega)) \times L^2(0, T;\hdiv)$ are solving Problem $P$, then $\Theta,p  \in L^2(0, T; H^1_0(\Omega))$ and $(\Theta,p)$ is a solution of Problem PC.
\end{proposition}

\beginproof. $ "\Rightarrow"$ Clearly, $\Theta$ and $p$ defined above have the regularity required in Problem $PC$. Furthermore, $\bu$ and $\bq$ are elements of $L^2(0, T; L^2(\Omega)^d)$. Recalling that $s(\cdot)$ is Lipschitz continuous, one immediately obtains that $s_C \in L^2(0, T; H^1_0(\Omega))$. Since $\partial_t s_C \in L^2(0, T; H^{-1}(\Omega))$ this shows that $s_C \in C(0, T; L^2(\Omega))$. With $t \in (0, T]$ and $\phi \in H^1_0(\Omega)$ arbitrary chosen, taking now $v = \chi_{(0, t]} \phi$ in \eqref{eq_cf_cont_1} and using the definition of $\bq$ gives
\begin{equation}\label{eq:prop3.1}
\la s(\Theta_C)-s(\Theta^0), \phi \ra - \la \int_0^t \bq(y) \, dy, \nabla \phi\ra = 0,
\end{equation}
for all $\phi \in H^1_0(\Omega)$. In other words, $\nabla \cdot \int_0^t \bq(y) \, dy = s(\Theta^0) - s(\Theta_C(t))$ in distributional sense. The regularity of $s_C$ mentioned above implies that, actually, $\int_0^t \bq (y) dy$ lies in $H^1(0, T; L^2(\Omega)) \cap L^2(0, T; H^1(\Omega))$ as well, and that $\int_0^t \bq  \in C(0, T; \hdiv)$. Moreover, by density arguments \eqref{eq_cont_1} holds for any $w \in L^2(\Omega)$. Also, from \eqref{eq:prop3.1} one gets $s_C \in C(0, T; L^2(\Omega))$. 

In a similar way, using \eqref{eq_cf_cont_3} and the definition of $\bu$ one obtains
$$
\big \la -\bu, \nabla w \big\ra = \la f_2(s_C), w \ra ,
$$
for all $w \in H^1_0(\Omega)$, so $\nabla \cdot \bu = f_2(s_C)$ for a.~e. $t \in [0, T]$. In view of the regularity of $s_C$ and of the assumptions on $f_2$, this means that $\bu \in C(0, T; \hdiv)$ and that \eqref{eq_cont_3} holds for every $t \in [0, T]$.

To obtain \eqref{eq_cont_4} one uses \eqref{eq_cf_cont_3}, the definition of $\bu$ and that $p_C$ has a vanishing trace on $\Gamma$. This gives
$$
\la a(s_C) \bu, \bv \ra =
- \la \nabla p_C + \bff_3(s_C), \bv \ra =
\la p_C , \nabla \cdot \bv \ra - \la \bff_3(s_C), \bv \ra ,
$$
for a.~e. $t \in [0, T]$ and for all $\bv \in \hdiv$. Recalling the continuity in time for $s_C$ and $\bu$, it follows that \eqref{eq_cont_4} is valid for every $t \in [0, T]$. Finally, one can use similar ideas to show that \eqref{eq_cont_2} holds true as well.\\
\mbox{ }\finishedproof $"\Rightarrow"$

$"\Leftarrow"$ Let now $(\Theta,\bq,p, \bu)$ be the solution of Problem $P$. We have to show that $\Theta,p \in L^2((0,T) \times \Omega)$ is the solution of $PC$, i.e. to show that the functions are actually in $ L^2(0,T; H_0^1(\Omega))$ and that they satisfy \eqref{eq_cf_cont_1}-\eqref{eq_cf_cont_3}. Further, since in (A1) $s(\cdot)$ is assumed Lipschitz, it follows that $s(\Theta) \in L^2(0, T; H^1(\Omega))$ as well. Clearly, \eqref{eq_cont_1} gives $\partial_t s(\Theta) = - \nabla \cdot \bq$ for a.e. $t \in [0, T]$. Since $\bq \in L^2(0, T; (L^2(\Omega))^d)$ one gets $\partial_t s(\Theta) \in L^2(0, T; H^{-1}(\Omega))$. 

Taking $v \in (C^\infty_0(\Omega))^d \subset H({\rm div},\Omega)$ arbitrary as test function in \eqref{eq_cont_4} we get
\begin{equation} \label{eq:proof:equivalence:1}
\la a(s(\Theta(y))) \bu, \bv \ra + \la \nabla p,  \bv \ra +  \la \bff_3(s(\Theta(y))), \bv \ra = 0,
\end{equation}
which implies
\begin{equation}\label{eq:proof:equivalence:2}
 \nabla p = - a(s(\Theta(y))) \bu - \bff_3(s(\Theta(y)))
\end{equation}
first in a distributional sense, and then by using the regularity of $\bu $ and $s(\cdot)$ in $L^2$ sense. It follows that $p \in L^2(0,T;H^1(\Omega))$. It remains to verify that $p$ has a vanishing trace on the boundary of $\Omega$. Taking now $\bv \in \hdiv$ and using the regularity of $p$ and \eqref{eq_cont_4} one obtains
\begin{equation} \label{eq:proof:equivalence:3}
 \la \nabla p,  \bv \ra  = - \la a(s(\Theta(y))) \bu, \bv \ra - \la \bff_3(s(\Theta(y))), \bv \ra =  - \la p, \nabla \cdot \bv \ra ,
\end{equation}
for all $v \in \hdiv$. By using now the Green theorem for $\bv \in \hdiv$, see \cite{brezzi}, pg. 91
$$
\int_{\Gamma} p \bv \cdot {\bf n} ds = (\nabla p, \bv) + (p, \nabla \cdot \bv) = 0.
$$
It follows immediately that $p \in  L^2(0,T; H^1_0(\Omega))$. From \eqref{eq:proof:equivalence:2} and \eqref{eq_cont_3} one gets that $p$ satisfies \eqref{eq_cf_cont_3}. In a similar manner one can show that $\Theta \in  L^2(0,T; H^1_0(\Omega))$ and it satisfies \eqref{eq_cf_cont_1}.\\
\mbox{ }\finishedproof $"\Leftarrow"$

We can now state an analogous result for the semi-discrete case.

{\bf Problem $PC^n$: Semi-discrete conformal variational formulation.} Let  $n \in \N, n \ge 1$ and  $\Theta^{n-1}$ be given. Find $\Theta_C^n$, $p_C^n  \in H^1_0(\Omega)$, such that  for all $v, w \in H^1_0(\Omega)$ one has
\begin{eqnarray}
\la s(\Theta_C^n) - s(\Theta_C^{n-1}), v \ra + \la \nabla \Theta_C^n + \frac {f_w(s\Theta_C^n)}{a(s(\Theta_C^n))}\nabla p_C^n, \nabla v\big\ra && \nonumber \\
+ \la \frac {f_w(s\Theta_C^n)}{a(s(\Theta_C^n))}\bff_3(s(\Theta_C^n)) - \bff_1(s(\Theta_C^n)), \nabla v\big\ra &=& 0, \label{eq_cf_semi_1}\\[1ex]\big \la \frac 1 {a(s(\Theta_C^n))} \big(\nabla p_C^n + \bff_3(s(\Theta_C^n))\big), \nabla w \big\ra &=& -\la f_2(s(\Theta_C^n)), w \ra. \label{eq_cf_semi_3}
\end{eqnarray}

\begin{proposition}\label{proposition:equivalence:semi}
Let $n \in \N, n \ge 1$, $\Theta^{n-1}$ be given and assume that (A1)-(A5) hold true. If $\Theta^n, p^n  \in  H^1_0(\Omega)$ is a solution to Problem $PC^n$ then, a solution to Problem $P^n$ is given by
$\Theta^n = \Theta_C^n$, $p^n = p_C^n$, $\bq^n = -\nabla \Theta_C^n -  \frac {f_w (s_C^n)} {a(s_C^n)}\left( \nabla p_C^n + \bff_3(s_C^n) \right) + \bff_1(s_C^n)$ and $\bu^n = - \frac 1 {a(s_C^n)}\left( \nabla p_C^n + \bff_3(s_C^n) \right)$. Conversely, if $(\Theta^n,\bq^n)  \in L^2(\Omega) \times \hdiv$,  $(p, \bu) \in L^2(\Omega) \times \hdiv$ are solving Problem $P^n$, then $\Theta^n,p^n  \in H^1_0(\Omega)$ and $(\Theta^n,p^n)$ is a solution of Problem $PC^n$.
\end{proposition}

The proof of Proposition \ref{proposition:equivalence:semi} is similar with the one of Proposition \ref{proposition:equivalence} (also see \cite{radu_SIAM}, Prop. 2.3) and it will be skipped here. Then the existence and uniqueness of a solution for the semi-discrete variational formulation \eqref{eq_semi-discrete_1}--\eqref{eq_semi-discrete_4} is a direct consequence of the equivalence result.

The following proposition establish the uniqueness for the fully-discrete variational problem \eqref{eq1_discrete_nonlin}-\eqref{eq4_discrete_nonlin}. As explained in the introduction of Sec. \ref{sec:equivalence}, the existence of a solution will follow from the convergence of the linear iteration scheme.

\begin{proposition}\label{prop:existence:fullydiscrete}
Let $n \in \N, n \ge 1$ be fixed. If (A1)-(A5) hold true and the time step $\tau$ is sufficiently small, then the problem \eqref{eq1_discrete_nonlin}-\eqref{eq4_discrete_nonlin} has at most one solution.
\end{proposition}

\beginproof. Let us assume that there exists two solutions of  \eqref{eq1_discrete_nonlin}-\eqref{eq4_discrete_nonlin}, $(\Theta^{n}_{h,i},\bq^{n}_{h,i},p^{n}_{h,i},\bu^{n}_{h,i}) \in W_h\times V_h \times W_h\times V_h$ with $i =1,2$. Further, we let $s^{n}_{h,i}:= s(\Theta^{n}_{h,i})$ stand for the two saturations.  These solutions are then satisfying
\begin{eqnarray}
\la s^{n}_{h,i}- s^{n-1}_h, w_h \ra + \tau \la \nabla \cdot \bq^{n}_{h,i}, w_h \ra &=& 0, \label{proof:ex:eq:1}\\[1ex]
\la \bq^{n}_{h,i}, \bv_h \ra -\la \Theta^{n}_{h,i},  \nabla  \cdot \bv_h \ra - \la f_w (s^{n}_{h,i}) \bu^{n}_{h,i},  \bv_h \ra  &=& \la \bff_1 (s^{n}_{h,i}), \bv_h \ra,\label{proof:ex:eq:2} \\[1ex]
\la \nabla \cdot \bu^{n}_{h,i}, w_h \ra   &=& \la f_2(s^{n}_{h,i}), w_h \ra, \label{proof:ex:eq:3}\\[1ex]
\la a(s^{n}_{h,i}) \bu^{n}_{h,i}, \bv_h \ra - \la p^{n}_{h,i}, \nabla \cdot \bv_h \ra + \la \bff_3(s^{n}_{h,i}), \bv_h \ra&=& 0\label{proof:ex:eq:4}
\end{eqnarray}
for $i =1,2$ and for all $w_h \in W_h$, $\bv_h \in V_h$. By subtracting \eqref{proof:ex:eq:3} and  \eqref{proof:ex:eq:4} for  $i = 2$ from the same equations for $i = 1$ we get for all $w_h \in W_h$, $\bv_h \in V_h$
\begin{eqnarray}
\la \nabla \cdot (\bu^{n}_{h,2}-\bu^{n}_{h,1}), w_h \ra   &=& \la f_2(s^{n}_{h,2}) -  f_2(s^{n}_{h,1}) , w_h \ra, \label{proof:ex:eq:5}\\[1ex]
\la a(s^{n}_{h,2}) \bu^{n}_{h,2} - a(s^{n}_{h,1}) \bu^{n}_{h,1},  \bv_h \ra - \la p^{n}_{h,2} - p^{n}_{h,1}, \nabla \cdot \bv_h \ra &=& \la \bff_3(s^{n}_{h,1})-  \bff_3(s^{n}_{h,2}),\bv_h\ra. \label{proof:ex:eq:6}
\end{eqnarray}
We test now \eqref{proof:ex:eq:5} with $w_h = p^{n}_{h,2} - p^{n}_{h,1} \in W_h$ and  \eqref{proof:ex:eq:6} with $\bv_h = \bu^{n}_{h,2} -  \bu^{n}_{h,1} \in V_h$, and add the resulting equations to obtain
\begin{equation}\label{proof:ex:eq:7}
\begin{array}{r}
\ds \la a(s^{n}_{h,2}) \bu^{n}_{h,2} - a(s^{n}_{h,1}) \bu^{n}_{h,1}, \bu^{n}_{h,2} -  \bu^{n}_{h,1}   \ra =  \la f_2(s^{n}_{h,2}) -  f_2(s^{n}_{h,1}) , p^{n}_{h,2} - p^{n}_{h,1} \ra \\[1ex]
\ds \hspace*{2cm}+  \la \bff_3(s^{n}_{h,1})-  \bff_3(s^{n}_{h,2}),\bu^{n}_{h,2} -  \bu^{n}_{h,1} \ra.
\end{array}
\end{equation}
After some algebraic manipulations, using (A2)-(A4) and Cauchy-Schwarz and Young inequalities one gets from \eqref{proof:ex:eq:7} above
\begin{equation}\label{proof:ex:eq:8}
\begin{array}{r}
\dfrac{a_\star}{4} \| \bu^{n}_{h,2} -  \bu^{n}_{h,1}  \|^2 \le (\dfrac{L_{f_2}^2}{ 2 \delta}  + \dfrac{L_{ \bff_3}^2}{ 2 a_\star} + \dfrac{L_{ a}^2M_u^2}{ a_\star})   \| s^{n}_{h,2} -  s^{n}_{h,1} \|^2 + \dfrac{ \delta}{ 2} \| p^{n}_{h,2} - p^{n}_{h,1} \|^2,
\end{array}
\end{equation}
for all $\delta > 0$. Using Lemma \ref{lemma_thomas} there exists $\bv_h \in V_h$ such that $\nabla \cdot \bv_h = p^{n}_{h,2} - p^{n}_{h,1} $ and $\|  \bv_h \| \le C_{\Omega,d}  \| p^{n}_{h,2} - p^{n}_{h,1} \|$. Testing with this $\bv_h $ in \eqref{proof:ex:eq:6} one gets
\begin{equation}\label{proof:ex:eq:9}
\begin{array}{r}
\| p^{n}_{h,2} - p^{n}_{h,1} \|^2 = \la a(s^{n}_{h,2}) \bu^{n}_{h,2} - a(s^{n}_{h,1}) \bu^{n}_{h,1},  \bv_h \ra + \la \bff_3(s^{n}_{h,1})-  \bff_3(s^{n}_{h,2}),\bv_h \ra.
\end{array}
\end{equation}
Using the properties of $\bv_h $, (A2)-(A4) and Cauchy-Schwarz's inequality,  \eqref{proof:ex:eq:9} implies
\begin{equation}\label{proof:ex:eq:10}
\| p^{n}_{h,2} - p^{n}_{h,1} \| \le C (\| \bu^{n}_{h,2} - \bu^{n}_{h,1} \| + \| s^{n}_{h,1}-  s^{n}_{h,2} \|),
\end{equation}
with $C$ not depending on the solutions or discretization parameters. From \eqref{proof:ex:eq:8} and  \eqref{proof:ex:eq:10} follows, by properly chosen $\delta$
\begin{equation}\label{proof:ex:eq:11}
\| \bu^{n}_{h,2} -  \bu^{n}_{h,1} \| \le C \| s^{n}_{h,1}-  s^{n}_{h,2} \|,
\end{equation}
with $C$ not depending on the solutions or discretization parameters. We proceed by subtracting \eqref{proof:ex:eq:1} and  \eqref{proof:ex:eq:2} for  $i = 2$ from the same equations for $i = 1$ to obtain
\begin{eqnarray}
\la s^{n}_{h,2} -s^{n}_{h,1} , w_h \ra + \tau \la \nabla \cdot (\bq^{n}_{h,2} - \bq^{n}_{h,1}), w_h \ra &=& 0, \label{proof:ex:eq:12}\\[1ex]
\la \bq^{n}_{h,2} -  \bq^{n}_{h,1}, \bv_h \ra -\la \Theta^{n}_{h,2} - \Theta^{n}_{h,1} ,  \nabla \cdot \bv_h \ra  &=& \la f_w (s^{n}_{h,2}) \bu^{n}_{h,2} -  f_w (s^{n}_{h,1}) \bu^{n}_{h,1},  \bv_h \ra \nonumber \\
&& + \la \bff_1 (s^{n}_{h,2}) - \bff_1 (s^{n}_{h,1}), \bv_h \ra,\label{proof:ex:eq:13}
\end{eqnarray}
for all $w_h \in W_h, \bv_h \in V_h$. Testing \eqref{proof:ex:eq:12} with $w_h = \Theta^{n}_{h,2} - \Theta^{n}_{h,1} \in W_h $ and  \eqref{proof:ex:eq:13} with $\bv_h =  \tau (\bq^{n}_{h,2} -  \bq^{n}_{h,1}) \in V_h $ and then adding the results gives
\begin{equation}\label{proof:ex:eq:14}
\begin{array}{l}
\ds \la s^{n}_{h,2} -s^{n}_{h,1} , \Theta^{n}_{h,2} - \Theta^{n}_{h,1}\ra + \tau \| \bq^{n}_{h,2} -  \bq^{n}_{h,1} \|^2 = \\[1ex]
\ds \hspace*{0.5cm}  \tau \la f_w (s^{n}_{h,2}) \bu^{n}_{h,2} -  f_w (s^{n}_{h,1}) \bu^{n}_{h,1},  \bq^{n}_{h,2} -  \bq^{n}_{h,1}\ra + \tau \la \bff_1 (s^{n}_{h,2}) - \bff_1 (s^{n}_{h,1}),  \bq^{n}_{h,2} -  \bq^{n}_{h,1} \ra.
\end{array}
\end{equation}
By some algebraic manipulations, using (A1)-(A4), the Cauchy-Schwarz and Young inequalities and the result \eqref{proof:ex:eq:11}, we obtain from \eqref{proof:ex:eq:14} above
\begin{equation}\label{proof:ex:eq:15}
\la s^{n}_{h,2} -s^{n}_{h,1} , \Theta^{n}_{h,2} - \Theta^{n}_{h,1}\ra + \dfrac{\tau}{4} \| \bq^{n}_{h,2} -  \bq^{n}_{h,1} \|^2 \le C \tau \|s^{n}_{h,1}-  s^{n}_{h,2}  \|^2,
\end{equation}
with $C$ not depending on the solutions or discretization parameters. Due to (A1), there holds $\ds \la s^{n}_{h,2} - s^{n}_{h,1} , \Theta^{n}_{h,2} - \Theta^{n}_{h,1}\ra \ge \dfrac{1}{L_s} \|s^{n}_{h,1}-  s^{n}_{h,2}  \|^2$, which together with \eqref{proof:ex:eq:15} immediately implies (for sufficiently small $\tau$) that $\bq^{n}_{h,2}  =  \bq^{n}_{h,1}$ and $ s^{n}_{h,1} = s^{n}_{h,2} $. Using these and \eqref{proof:ex:eq:11}, we obtain that $\bu^{n}_{h,2}  =  \bu^{n}_{h,1}$. From \eqref{proof:ex:eq:10} results also $p^{n}_{h,2} = p^{n}_{h,1} $. Finally, equation \eqref{proof:ex:eq:13} will furnish $\Theta^{n}_{h,2} =\Theta^{n}_{h,1}$, which completes the proof of uniqueness.

\mbox{ }\finishedproof

\begin{remark} The uniqueness of a solution for the $L$-scheme introduced in \eqref{eq1_discrete_lin}--\eqref{eq4_discrete_lin} can be proved similarly. Since this is a linear algebraic system, uniqueness also gives the existence of a solution.
\end{remark}

\subsection{Stability estimates}
As in \cite{radu_NumMath} and, actually, as in Proposition \ref{proposition:equivalence:semi} one can obtain some stability estimates for the Problems $P^n$. Moreover, the same holds for Problem $P^n_h$, but these estimates are not needed for proving the convergence of the scheme and are therefore skipped here. Following  \cite{radu_NumMath}, Lemma 3.2, pg. 293 there holds:
\begin{proposition}\label{prop:stability}
Assume that (A1)-(A5) hold true. Let $(\Theta^n, \bq^n, p^n, \bu^n)$, $n \in \N, n \ge 1$ be the solution of Problem $P^n$. Then there holds
\begin{eqnarray}
 \tau \sum_{n = 1}^N  \| \Theta^n \|_1^2 + \tau \sum_{n = 1}^N  \|  \bq^n \|^2_{div} + \tau \sum_{n = 1}^N  \| p^n \|_1^2 + \tau \sum_{n = 1}^N  \|  \bu^n \|^2_{div}&\le& C,
 \end{eqnarray}
with $C$ not depending on discretization parameters.
\end{proposition}

\subsection{A priori error estimates}
Having established the existence and uniqueness for Problems $P, P^n$ and $P_h^n$, we can focus now on the convergence of the scheme  (\ref{eq1_discrete_nonlin})-(\ref{eq4_discrete_nonlin}). This will be done by deriving {\it a priori} error estimates. We assume that the fully discrete non-linear problem (\ref{eq1_discrete_nonlin})-(\ref{eq4_discrete_nonlin}) is solved exactly. The proofs of this section follow the lines in \cite{radu_NumMath} and \cite{chen2}. The following two propositions will quantify the error between the continuous and the semi-discrete formulations, and between the semi-discrete and discrete ones, respectively. Finally the two propositions will be put together to obtain  the main convergence result given in Theorem \ref{th_main_result_estimates}.

Recalling the definition $\bar{g}^n := \frac 1 \tau \int_{t_{n-1}}^{t_n} g(t) dt \in X,$ for any $g \in L^2(0, T; X)$ and $X$ a Banach space, we have
\begin{proposition} \label{prop1} Let $(\Theta,\bq, p, \bu)$ be the solution of Problem $P$ and $(\Theta^n, \bq^n, p^n, \bu^n)$ be the solution of $P^n, n \in \N, n\ge 1.$ Assuming (A1)-(A5) and that the time step is small enough there holds:
\begin{displaymath}
\sum_{n=1}^N \int_{t_{n-1}}^{t_{n}} \la s(\Theta(t)) - s(\Theta^n), \Theta(t) - \Theta^n \ra \, dt +  \| \sum_{n=1}^N \int_{t_{n-1}}^{t_{n}} \bq - \bq^n \, dt\|^2
\end{displaymath}
\vspace*{-0.5cm}
\begin{eqnarray}
\hspace*{2cm} + \sum_{n=1}^N  \| \int_{t_{n-1}}^{t_{n}} (\bq - \bq^n) \, dt \|^2  & \le & C \tau, \label{prop1_result_1}\\[2ex]
\hspace*{2cm}  \tau \| \overline{\bu}^n - \bu^n \|^2 + \tau \| \overline{p}^n - p^n \|^2  & \le &  C \int_{t_{n-1}}^{t_{n}}  \| s(\Theta(t)) - s(\Theta^n)\|^2 dt, \label{prop1_result_2}\\[2ex]
 \| \sum_{n=1}^N \int_{t_{n-1}}^{t_{n}} (\Theta(t) - \Theta^n) \, dt\|^2 &\le& C \tau, \label{prop1_result_3}
\end{eqnarray}
with the constants $C$ not depending on the discretization parameters.
\end{proposition}

\beginproof. We start with proving \eqref{prop1_result_2}. By integrating \eqref{eq_cont_3}, \eqref{eq_cont_4} from $t_{n-1}$ to $t_n$ one obtains
\begin{eqnarray}
\la \nabla \cdot \bun, w \ra &=& \la \overline{f_2(s)}^n, w \ra, \label{prop1_eq:1}\\[1ex]
\la \overline{a(s) \bu}^n, \bv \ra - \la \opn, \nabla \cdot \bv \ra +  \la \overline{\bff_3(s)}^n, \bv \ra &=& 0, \label{prop1_eq:2}
\end{eqnarray}
for all $w \in L^2(\Omega)$ and $\bv \in \hdiv$. By subtracting now  \eqref{eq_semi-discrete_3} and \eqref{eq_semi-discrete_4} from  \eqref{prop1_eq:1} and \eqref{prop1_eq:2}, respectively, we get
\begin{eqnarray}
\la \nabla \cdot (\bun - \bu^n), w \ra &=& \la \overline{f_2(s)}^n - f_2(s^n), w \ra, \label{prop1_eq:1a}\\[1ex]
\la \overline{a(s) \bu}^n - a(s^n) \bu^n, \bv \ra - \la \opn - p^n, \nabla \cdot \bv \ra &=& - \la \overline{\bff_3(s)}^n -\bff_3(s^n), \bv \ra, \label{prop1_eq:2a}
\end{eqnarray}
for all $w \in L^2(\Omega)$ and $\bv \in \hdiv$. Taking $w = \opn - p^n \in L^2(\Omega)$ in \eqref{prop1_eq:1a} and $\bv = \bun - \bu^n \in \hdiv$ in \eqref{prop1_eq:2a}, and summing the results we obtain
\begin{equation}\label{prop1_eq:3}
\la \overline{a(s) \bu}^n - a(s^n)\bu^n, \bun - \bu^n \ra = \la \overline{f_2(s)}^n - f_2(s^n), \opn - p^n \ra - \la \overline{\bff_3(s)}^n - \bff_3(s^n),  \bun - \bu^n \ra.
\end{equation}
By Young's inequality, this further implies
\begin{displaymath}
\begin{array}{l}
\ds \la  \dfrac{1}{\tau} \int_{t^{n-1}}^{t_n} a(s) \bu -  a(s^n) \bu^n \, dt, \bun - \bu^n \ra \le \dfrac{1}{2\delta} \|  \overline{f_2(s)}^n - f_2(s^n) \|^2 + \dfrac{\delta}{2} \| \opn - p^n \|^2  \\[1ex]
\ds \hspace*{4cm}+ \dfrac{1}{a_\star} \| \overline{\bff_3(s)}^n - \bff_3(s^n) \|^2 + \dfrac{a_\star}{4} \| \bun - \bu^n \|^2.
\end{array}
\end{displaymath}
The above is further equivalent to 
$$
\begin{array}{l}
\ds \la \dfrac{1}{\tau} \int_{t^{n-1}}^{t_n} (a(s) -  a(s^n)) \bu \, dt, \bun - \bu^n \ra  +  \la \dfrac{a(s^n) }{\tau} \int_{t^{n-1}}^{t_n} \bu -  \bu^n \, dt, \bun - \bu^n \ra \\[1ex]
\hspace*{2cm}\le \dfrac{1}{2\delta} \|  \overline{f_2(s)}^n - f_2(s^n) \|^2
+ \dfrac{\delta}{2} \| \opn - p^n \| ^2+ \dfrac{1}{a_\star} \| \overline{\bff_3(s)}^n - \bff_3(s^n) \|^2 + \dfrac{a_\star}{4} \| \bun - \bu^n \|^2,
\end{array}
$$
which by using (A2)-(A3) leads to
\begin{equation}\label{prop1_eq:4}
\ds 
\dfrac{3a_{\star}}{4} \| \bun - \bu^n \|^2 \le  (\dfrac{L_{f_2}^2 }{2 \tau \delta} + \dfrac{L_{\bff_3}^2 }{2 \tau a_\star}) \int_{t_{n-1}}^{t_n} \|  s - s^n \|^2 \, dt  + \dfrac{\delta}{2} \| \opn - p^n \| ^2 - T_1,
\end{equation}
for all $\delta > 0$, where $\ds T_1 = \la \dfrac{1}{\tau} \int_{t^{n-1}}^{t_n} (a(s) -  a(s^n)) \bu \, dt, \bun - \bu^n \ra$. Now, to estimate $T_1$ one uses (A2), (A4) and the Young inequality to obtain
\begin{eqnarray}
| T_1 | &\le & \| \dfrac{1}{\tau} \int_{t^{n-1}}^{t_n} (a(s) -  a(s^n)) \bu \, dt \| \| \bun - \bu^n \| \nonumber \\[1ex]
&\le &\dfrac{1}{\tau²^2 a_\star} \int_\Omega \left( \int_{t^{n-1}}^{t_n} (a(s) -  a(s^n)) \bu \, dt \right)^2 \, dx + \dfrac{a_\star}{4}  \| \bun - \bu^n \|^2 \nonumber \\[1ex]
&\le & \dfrac{M^2_\bu}{\tau² a_\star} \int_\Omega \int_{t^{n-1}}^{t_n} \left( (a(s) -  a(s^n) \right)^2 \, dt \, dx + \dfrac{a_\star}{4}  \| \bun - \bu^n \|^2  \nonumber \\[1ex]
&\le &  \dfrac{M^2_\bu L^2_a}{\tau² a_\star} \int_{t^{n-1}}^{t_n} \| s -  s^n \|^2 \, dt + \dfrac{a_\star}{4}  \| \bun - \bu^n \| ^2. \label{prop1_eq:5}
\end{eqnarray}
From \eqref{prop1_eq:4} and \eqref{prop1_eq:5} it follows immediately
\begin{equation}\label{prop1_eq:6}
 \dfrac{a_\star}{2} \| \overline{\bu}^n - \bu^n \|^2  \le \dfrac{C}{\tau} \int_{t_{n-1}}^{t_{n}}  \| s(\Theta(t)) - s(\Theta^n) \|^2 dt + \dfrac{\delta}{2} \| \opn - p^n \|^2,
\end{equation}
with $C$ not depending on the discretization parameters.

To estimate the last term above, one uses Lemma \ref{lemma_thomas_cont}, ensuring the existence of a $\bv \in \hdiv $ such that $\nabla \cdot \bv = \opn - p^n$  and $\| \bv \| \le C_\Omega \| \opn - p^n \| $. Using this as test function in \eqref{prop1_eq:2a} gives
\begin{eqnarray}
\| \opn - p^n \|^2 &=&  \la  \overline{a(s) \bu}^n - a(s^n)\bu^n, \bv \ra +  \la \overline{\bff_3(s)}^n -\bff_3(s^n), \bv \ra \nonumber \\[1ex]
& \leq & C^2_\Omega \| \overline{a(s) \bu}^n - a(s^n)\bu^n \|^2 + C^2_\Omega \|  \overline{\bff_3(s)}^n -\bff_3(s^n) \|²^2 + \dfrac{1}{2}\| \opn - p^n \|^2. \label{prop1_eq:7}
\end{eqnarray}
Proceeding as for \eqref{prop1_eq:6}, for \eqref{prop1_eq:7} we get:
\begin{equation}\label{prop1_eq:8}
\| \opn - p^n \|^2 \le \dfrac{C}{\tau} \int_{t_{n-1}}^{t_{n}}  \| s(\Theta(t)) - s(\Theta^n) \|^2 dt + C \|  \overline{\bu}^n - \bu^n \|^2,
\end{equation}
with the constant $C$ not depending on the discretization parameters. Using \eqref{prop1_eq:8} and \eqref{prop1_eq:6}, and choosing $\delta$ properly, one obtains \eqref{prop1_result_2}.

To prove \eqref{prop1_result_1} we follow the steps in the proof of Lemma 3.3, pg. 296 in \cite{radu_NumMath}. By summing up \eqref{eq_semi-discrete_1} for $k = 1$ to $n$ and subtracting \eqref{eq_cont_1} from the resulting we get for all $w \in \l2$
\begin{equation} \label{prop1_eq:9}
\la s(\Theta(t_n)) - s^n, w \ra + \tau \sum_{k=1}^n \la \nabla \cdot (\oqn - \bq^k), w \ra = 0.
\end{equation}
Further,  subtracting \eqref{eq_cont_2} at $t = t_{k-1}$ from \eqref{eq_cont_2} at $t = t_{k}$, dividing by the time step size $\tau$ and subtracting from the result  \eqref{eq_semi-discrete_2}  we obtain for all $\bv \in \hdiv$
\begin{equation}
\la \oqn - \bq^n, \bv \ra - \la \othetan - \Theta^n, \nabla \cdot \bv \ra - \la \overline{f_w(s) \bu}^n - f_w(s^n) \bu^n, \bv \ra  =  \la \overline{\bff_1(s)}^n -\bff_1(s^n), \bv \ra. \label{prop1_eq:10}
\end{equation}
By testing \eqref{prop1_eq:9} with $w =  \othetan - \Theta^n \in \l2$ and \eqref{prop1_eq:10} with $\ds \bv = \tau \sum_{k=1}^n (\oqk - \bq^k) \in \hdiv$, adding the results and summing up from $n = 1$ to $N$ we get
\begin{equation}\label{prop1_eq:11}
\begin{array}{l}
 \ds \sum_{n=1}^N \la s(\Theta(t_n)) - s^n,  \othetan - \Theta^n \ra + \sum_{n=1}^N \tau \la \oqn - \bq^n, \sum_{k=1}^n (\oqk - \bq^k) \ra \\[1ex]
\hspace*{2cm}\ds  -  \sum_{n=1}^N \la  \overline{f_w(s) \bu}^n - f_w(s^n) \bu^n, \tau \sum_{k=1}^n (\oqk - \bq^k) \ra =  \sum_{n=1}^N  \la \overline{\bff_1(s)}^n -\bff_1(s^n), \tau \sum_{k=1}^n (\oqk - \bq^k) \ra .
\end{array}
\end{equation}
We estimate separately each term in \eqref{prop1_eq:11}, which are denoted by $T_1$, $T_2$, $T_3$ and  $T_4$. For $T_1$ we proceed as in the proof of Lemma 3.3, pg. 296 in \cite{radu_NumMath}, as the term here is identical to the one there and obtain
\begin{equation}\label{prop1_eq:12}
T_1 = \dfrac{1}{\tau} \sum_{n=1}^N \int_{t_{n-1}}^{t_n} \la s(\Theta) - s^n,  \Theta - \Theta^n \ra \, dt + \dfrac{1}{\tau} \sum_{n=1}^N\int_{t_{n-1}}^{t_n} \la s (\Theta(t_n))  - s(\Theta), \Theta - \Theta^n \ra \, dt
\end{equation}
with the first term above being positive to remain on the left hand side of \eqref{prop1_result_1}. Using the regularity of the solutions (both continuous and semi-discrete) and the stability estimates in Proposition \ref{prop:stability}) one can follow the steps in estimating $T_{11}$ in \cite{radu_NumMath} to obtain for the second term above
\begin{equation}\label{prop1_eq:13}
 | \dfrac{1}{\tau} \sum_{n=1}^N\int_{t_{n-1}}^{t_n} \la s (\Theta(t_n))  - s(\Theta), \Theta - \Theta^n \ra | \le C,
\end{equation}
with $C$ not depending on the discretization parameters. Moreover, if the data is such that both phases present at any time and everywhere in the system, the estimate in \eqref{prop1_eq:13} can be improved to $C \tau$, as discussed in Remark \ref{corollary_optimal_estimates} below. Such estimates are optimal. 

For the second term in \eqref{prop1_eq:11} one uses the algebraic identity \eqref{lemma1_1} to obtain 

\begin{equation}\label{prop1_eq:14}
 T_2 = \dfrac{\tau}{2}\| \sum_{n=1}^N (\oqn - \bq^n) \|^2 +   \sum_{n=1}^N \dfrac{\tau}{2} \| \oqn - \bq^n\|^2.
\end{equation}
The two terms above will remain on the left hand side of \eqref{prop1_result_1}. We proceed by estimating $T_3$ in \eqref{prop1_eq:11}. By the Young inequality there holds
\begin{eqnarray}
| T_3 | &=& | \sum_{n=1}^N \la  \overline{f_w(s) \bu}^n - f_w(s^n) \bu^n, \tau \sum_{k=1}^n (\oqk - \bq^k) \ra | \nonumber \\[1ex]
& \le & \dfrac{\delta}{2} \sum_{n=1}^N \| \overline{f_w(s) \bu}^n - f_w(s^n) \bu^n \|^2 + \dfrac{\tau^2}{2\delta }\sum_{n=1}^N  \|   \sum_{k=1}^n (\oqk - \bq^k) \|²^2 \nonumber \\[1ex]
& \le & T_{31} + \dfrac{\tau^2}{2\delta }\sum_{n=1}^N  \|   \sum_{k=1}^n (\oqk - \bq^k) \|²^2.\label{prop1_eq:15}
\end{eqnarray}
The second term on the right is estimated by using the Gronwall lemma. For the first one, one uses (A3), (A4) and \eqref{prop1_result_2} to obtain
\begin{eqnarray}
| T_{31} | &=& \dfrac{\delta}{2} \sum_{n=1}^N \int_\Omega \left( \dfrac{1}{\tau} \int_{t_{n-1}}^{t_n} f_w(s) \bu - f_w(s^n) \bu^n \, dt \right)^2 dx \nonumber \\[1ex]
& \le & \dfrac{\delta}{\tau²^2 } \sum_{n=1}^N  \int_\Omega  \left( \int_{t_{n-1}}^{t_n} (f_w(s)  - f_w(s^n) )\bu\, dt \right)^2 dx +  \dfrac{\delta}{\tau²^2 } \sum_{n=1}^N  \int_\Omega f^2_w(s^n)  \left(  \int_{t_{n-1}}^{t_n} (\bu - \bu^n) \, dt \right)^2 dx  \nonumber \\[1ex]
& \le &  \dfrac{\delta M^2_\bu L^2_{f_w}}{\tau }  \sum_{n=1}^N \int_{t_{n-1}}^{t_n} \| s - s^n \|^2 +  \delta M^2_{f_w} \sum_{n=1}^N \| \bun - \bu^n \|^2  \nonumber \\[1ex]
& \le &  \dfrac{\delta M^2_\bu L^2_{f_w}}{\tau }  \sum_{n=1}^N \int_{t_{n-1}}^{t_n} \| s - s^n \|^2 +  \dfrac{\delta M^2_{f_w}C}{\tau} \sum_{n=1}^N \int_{t_{n-1}}^{t_n} \| s - s^n \|^2   \nonumber \\[1ex]
& \le & \dfrac{C \delta}{\tau} \sum_{n=1}^N \int_{t_{n-1}}^{t_n} \| s - s^n \|^2 \, dt,\label{prop1_eq:16}
\end{eqnarray}
for all $\delta > 0$ and with a constant $C$ not depending on the discretization parameters. In the same manner one can bound the last term in \eqref{prop1_eq:11}. Using again Young's inequality and (A5), one gets
\begin{equation}\label{prop1_eq:17}
| T_4 | \le \dfrac{C \delta^\prime}{\tau } \sum_{n= 1}^N \int_{t_{n-1}}^{t_n} \| s - s^n \|^2 \, dt  + \dfrac{\tau^2}{ 2 \delta^\prime} \sum_{n= 1}^N \| \sum_{k=1}^n (\oqk - \bq^k)\|^2
\end{equation}
for all $\delta^\prime > 0$ and with a constant $C$ not depending on the discretization parameters. We observe that due to (A1) there holds
\begin{equation}\label{prop1_eq:18}
 \sum_{n= 1}^N \int_{t_{n-1}}^{t_n} \| s - s^n \|^2 \le  \sum_{n= 1}^N \int_{t_{n-1}}^{t_n} \la s - s^n, \Theta - \Theta^n \ra \, dt.
\end{equation}
By using \eqref{prop1_eq:18}, choosing $\delta$ and $\delta^\prime$ properly, the first terms in the right hand sides of \eqref{prop1_eq:16} and \eqref{prop1_eq:17} are absorbed in the left hand side of \eqref{prop1_result_1}. Putting now together \eqref{prop1_eq:11} - \eqref{prop1_eq:17} and applying the discrete Gronwall lemma gives \eqref{prop1_result_1}.

Finally, to prove \eqref{prop1_result_3} one follows the step in Lemma 3.9, pg. 300 in \cite{radu_NumMath}. By Lemma \ref{lemma_thomas_cont}, there exists a function $ \bv \in \hdiv$ which satisfies $\ds \nabla \cdot \bv = \sum_{n=1}^N (\overline{{\Theta}}^n - \Theta^n)$ and $ \ds \| \bv \| \le C_\Omega \| \sum_{n=1}^N (\overline{{\Theta}}^n - \Theta^n)\|$. We use this as test function in \eqref{prop1_eq:10}. Now \eqref{prop1_result_3} follows from \eqref{prop1_result_1}.

\finishedproof\\

The next proposition quantifies the error between the semi-discrete solution and the fully discrete one. Recall the notations $s^k = s(\Theta^k)$ and $s_h^k = s(\Theta_h^k)$, $k \in \N$.
\begin{proposition}\label{prop2} Let $n \in \N, n\ge 1$ and let $(\Theta^n, \bq^n, p^n, \bu^n)$ be the solution of $P^n$, and $(\Theta^n_h, \bq^n_h, p^n_h, \bu^n_h)$ be the solution of $P^n_h$.  Assuming (A1)-(A5) and that the time step is small enough, there holds
\begin{equation}\label{prop2_result_1}
\begin{array}{l}
\ds \sum_{n=1}^N \left( \la s^n - s^n_h, \Theta^n - \Theta^n_h \ra + \| s^n - s^n_h\|^2 \right) + \tau \|  \sum_{n=1}^N (\Pi_h  \bq^n - \bq^n_h) \|^2  \\[1ex]
\ds \hspace*{1cm} \le  C \sum_{n=1}^N (\| \bq^n - \Pi_h \bq^n \|^2 +  \| \Theta^n - P_h \Theta^n \|^2 + \| \bu^n - \Pi_h \bu^n \|^2 + \| p^n - P_h p^n \|^2)
\end{array}
\end{equation}
and
\begin{equation}\label{prop2_result_2}
\begin{array}{l}
\ds \| \bu^n -\bu^n_h \|^2 + \| \nabla \cdot (\bu^n -\bu^n_h) \|^2 + \| p^n - p^n_h \|^2 \\[2ex]
\ds \hspace*{1cm} \le C(\| \bu^n - \Pi_h \bu^n \|^2 + \| s^n - s^n_h\|^2 + \|p^n - P_h p^n \|^2),
\end{array}
\end{equation}
with the constants $C$ above not depending on the discretization parameters.
\end{proposition}
\beginproof. The proof of \eqref{prop2_result_2} can be found in \cite{chen2}, where a MFEM was applied for the discretization of the pressure equation, but the Galerkin FEM for the saturation equation. Therefore we give here only the proof of \eqref{prop2_result_1}. By subtracting \eqref{eq1_discrete_nonlin} and \eqref{eq2_discrete_nonlin} from \eqref{eq_semi-discrete_1} and \eqref{eq_semi-discrete_2}, summing up from $k = 1$ to $n$ and using the properties of the projectors, one gets
\begin{eqnarray}
\la s^n - s^n_h, w_h \ra + \tau \sum_{k=1}^n \la \nabla \cdot (\Pi_h \bq^k - \bq^k_h), w_h \ra &=& 0, \label{prop2_eq:1}\\ [1ex]
\la \bq^n - \bq^n_h,\bv_h\ra - \la P_h \Theta^n - \Theta^n_h,\nabla \cdot \bv_h \ra - \la f_w(s^n) \bu^n - f_w(s^n_h) \bu^n_h, \bv_h \ra &=& \la \bff_1(s^n) - \bff_1(s^n_h) , \bv_h \ra\nonumber\\
\label{prop2_eq:2}
\end{eqnarray}
for all $w_h \in W_h$ and $\bv_h \in V_h.$  Taking $w_h = P_h \Theta^n -  \Theta^n_h \in W_h$ and $\bv_h = \tau \sum_{k=1}^n  (\Pi_h \bq^k - \bq^k_h) \in V_h $ in \eqref{prop2_eq:1} and \eqref{prop2_eq:2}, respectively, adding the results and summing up from $n = 1$ to $N$ we obtain
\begin{equation}\label{prop2_eq:3}
\begin{array}{l}
\ds \sum_{n=1}^N \la s^n - s^n_h, P_h \Theta^n -  \Theta^n_h \ra + \tau \sum_{n=1}^N \la \bq^n - \bq^n_h,  \sum_{k=1}^n (\Pi_h \bq^k - \bq^k_h) \ra  \\[1ex]
\ds \hspace*{1cm} - \sum_{n=1}^N  \la f_w(s^n) \bu^n - f_w(s^n_h) \bu^n_h, \tau \sum_{k=1}^n  (\Pi_h \bq^k - \bq^k_h) \ra =  \sum_{n=1}^N  \la \bff_1(s^n) - \bff_1(s^n_h), \tau \sum_{k=1}^n  (\Pi_h \bq^k - \bq^k_h) \ra . \end{array}
\end{equation}
Denoting the terms above by $\hat T_1$, $\hat T_2$, $\hat T_3$ and $\hat T_4$, we proceed by estimating them separately. For $\hat T_1$ there holds
\begin{equation}\label{prop2:_eq:5}
\hat T_1  =  \sum_{n=1}^N  \la s^n - s^n_h, \Theta^n - \Theta^n_h \ra + \sum_{n=1}^N  \la s^n - s^n_h, P_h \Theta^n - \Theta^n \ra
\end{equation}
with the first part above being positive due to the monotonicity of $s(\cdot)$.  By Young's inequality, for the second term in \eqref{prop2:_eq:5} one gets
\begin{equation}\label{prop2:_eq:5a}
\sum_{n=1}^N  \la s^n - s^n_h, P_h \Theta^n - \Theta^n \ra  \le  \dfrac{\delta_1}{2} \sum_{n=1}^N  \| s^n - s^n_h \|^2 + \dfrac{1}{2\delta_1} \sum_{n=1}^N  \| P_h \Theta^n - \Theta^n \|^2,
\end{equation}
for all $\delta_1 > 0$. Note that due to (A1) there holds
\begin{equation}\label{prop2_eq:6}
\sum_{n=1}^N  \la s^n - s^n_h, \Theta^n - \Theta^n_h \ra \ge \sum_{n=1}^N  \dfrac{1}{L_s} \| s^n - s^n_h \|^2.
\end{equation}
After properly choosing $\delta_1$, the second term on the right in \eqref{prop2:_eq:5a} can be absorbed by $ \frac 1 2 \sum_{n=1}^N  \la s^n - s^n_h, \Theta^n - \Theta^n_h \ra $.

Using the algebraic identity \eqref{lemma1_1}, for $\hat T_2$ it holds 
\begin{eqnarray}
\hat T_2 & = & \tau \sum_{n=1}^N \la \bq^n - \Pi_h \bq^n,  \sum_{k=1}^n (\Pi_h \bq^k - \bq^k_h) \ra + \tau \sum_{n=1}^N \la \Pi_h \bq^n - \bq^n_h,  \sum_{k=1}^n (\Pi_h \bq^k - \bq^k_h) \ra  \nonumber \\[1ex]
& = & \hat T_{21} + \dfrac{\tau}{2} \| \sum_{n=1}^N (\Pi_h \bq^n - \bq^n_h) \|^2 + \dfrac{\tau}{2} \sum_{n=1}^N \| \Pi_h \bq^n - \bq^n_h \|^2.\label{prop2_eq:7}
\end{eqnarray}
The only term remaining to be estimated is $\hat T_{21}$. This is done by using Young's inequality
\begin{equation}
| \hat T_{21} | \le    \dfrac{1}{2} \sum_{n=1}^N \| \bq^n - \Pi_h \bq^n \|^2 + \dfrac{\tau^2}{2} \sum_{n=1}^N \| \sum_{k=1}^n (\Pi_h \bq^k - \bq^k_h) \|^2  \label{prop2_eq:8}.
\end{equation}
In estimating $\hat T_3$  we use (A2)-(A4), Young's inequality and \eqref{prop2_result_2}. There holds
\begin{eqnarray}
| \hat T_{3} |  & = &  | \sum_{n=1}^N  \la f_w(s^n) \bu^n - f_w(s^n_h) \bu^n_h, \tau \sum_{k=1}^n  \Pi_h \bq^k - \bq^k_h \ra | \nonumber \\[1ex]
& \le & \dfrac{\delta_3}{2} \sum_{n=1}^N  \|  f_w(s^n) \bu^n - f_w(s^n_h) \bu^n_h \|^2 +  \dfrac{\tau^2}{2\delta_3} \sum_{n=1}^N  \| \sum_{k=1}^n (\Pi_h \bq^k - \bq^k_h) \|^2  \nonumber \\[1ex]
& \le &  C \delta_3 \sum_{n=1}^N \| s^n - s^n_h \|^2 + \delta_3 M^2_{f_w} \sum_{n=1}^N \| \bu^n - \bu^n_h \|^2  +  \dfrac{\tau^2}{2\delta_3} \sum_{n=1}^N  \| \sum_{k=1}^n (\Pi_h \bq^k - \bq^k_h) \|^2  \nonumber \\[1ex]
& \le &  C \delta_3 \sum_{n=1}^N \| s^n - s^n_h \|^2 + C \sum_{n=1}^N \| \bu^n - \bu^n_h \|^2 +  \dfrac{\tau^2}{2\delta_3} \sum_{n=1}^N  \| \sum_{k=1}^n (\Pi_h \bq^k - \bq^k_h) \|^2 \label{prop2_eq:9}
\end{eqnarray}
for all $\delta_3 >0$ and with the constants $C$ not depending on the discretization parameters. In a similar manner, by using (A5) we can bound also the last term $\hat T_4$. There holds for all $\delta_4 >0$
\begin{eqnarray}
| \hat T_{4} |  & \le  &  C \delta_4 \sum_{n=1}^N   \| s^n - s^n_h \|^2  + \dfrac{\tau^2}{\delta_4} \sum_{n=1}^N \| \sum_{k=1}^n (\Pi_h \bq^k - \bq^k_h) \|^2. \label{prop2_eq:10}
\end{eqnarray}
Finally, putting together \eqref{prop2_eq:3} - \eqref{prop2_eq:10}, choosing $\delta_1-\delta_4$ properly, and using the discrete Gronwall lemma we obtain the result \eqref{prop2_result_1}.

\finishedproof

The main result below is a straightforward consequence of Proposition \ref{prop1} and Proposition \ref{prop2}, the properties of the projectors and the regularity of the solution.
\begin{theorem}\label{th_main_result_estimates} Let $(\Theta,\bq, p, \bu)$ be the solution of Problem $P$ and let $(\Theta^n_h, \bq^n_h, p^n_h, \bu^n_h)$ be the solution of $P^n_h$, $n \in \{1, \ldots, N\}$.  Assuming (A1)-(A5) and that the time step is small enough, there holds
\begin{eqnarray}
\sum_{n=1}^N \int_{t_{n-1}}^{t_{n}} \la s(\Theta(t)) - s(\Theta^n_h), \Theta(t) - \Theta^n_h \ra \, dt + \| \sum_{n=1}^N \int_{t_{n-1}}^{t_{n}} \bq - \bq^n_h \, dt\|^2 \le C (\tau + h^2), \\[2ex]
\sum_{n=1}^N \tau \| \bun -  \bu^n_h \| ^2 + \sum_{n=1}^N \tau \| \opn -  p^n_h \| ^2 \le C (\tau + h^2),
\end{eqnarray}
with the constant $C$ not depending on the discretization parameters. \end{theorem}

\begin{remark} The error estimates presented above can be extended to the case of $s(\cdot)$ being only H\"older continuous (instead Lipschitz continuous), see  \cite{radu2016hoelder} for the details.
\end{remark}

\begin{remark} \label{corollary_optimal_estimates} In the non-degenerate case, when the disappearance of phases is not allowed, one can obtain the optimal error estimates (similar to Corollary 3.6, pg. 299 in \cite{radu_NumMath})
\begin{eqnarray}
\sum_{n=1}^N \int_{t_{n-1}}^{t_{n}} \la s(\Theta(t)) - s(\Theta^n_h), \Theta(t) - \Theta^n_h \ra \, dt + \| \sum_{n=1}^N \int_{t_{n-1}}^{t_{n}} \bq - \bq^n_h \, dt\|^2 \le C (\tau^2 + h^2), \\[2ex]
\sum_{n=1}^N \tau \| \bun -  \bu^n_h \| ^2 + \sum_{n=1}^N \tau \| \opn -  p^n_h \| ^2 \le C (\tau^2 + h^2).
\end{eqnarray}
\end{remark}

\section{Linearization scheme}\label{sec:lin}
In this section we analyze the convergence of the (fully discrete) linearization scheme \eqref{eq1_discrete_lin}--\eqref{eq4_discrete_lin} proposed to solve the non-linear system  (\ref{eq1_discrete_nonlin}) --  (\ref{eq4_discrete_nonlin}). We show that the scheme is robust it converges linearly. The scheme does not involve any regularization step. As the scheme is used to solve the nonlinear systems in one time step, throughout this section $n \in N, n \ge 1$ is fixed. For the ease of the presentation we recall the scheme \eqref{eq1_discrete_lin}--\eqref{eq4_discrete_lin}:

{\bf Problem $P_h^{n,i}$: Linearization scheme ($L$-scheme).} Let $L \ge L_s$, $i \in \N, i \ge 1$ and let  $\Theta^{n, i-1}_h \in W_h$ be given. Find $\Theta^{n,i}_h,  p^{n,i}_h \in W_h$ and $\bq^{n,i}_h, \bu^{n,i}_h \in V_h$ such that
\begin{eqnarray*}
\la L (\Theta^{n,i}_h - \Theta^{n,i-1}_h) + s^{n,i-1}_h, w_h \ra  + \tau \la \nabla \cdot \bq^{n,i}_h, w_h \ra &=& \la s^{n-1}_h, w_h \ra,\\
\la \bq^{n,i}_h, \bv_h \ra -\la \Theta^{n,i}_h,  \nabla  \cdot \bv_h \ra - \la f_w (s^{n,i-1}_h) \bu^{n,i}_h,  \bv_h \ra &=&  \la \bff_1(s^{n,i-1}_h) , \bv_h \ra ,\\
\la \nabla \cdot \bu^{n,i}_h, w_h \ra   &=& \la f_2(s^{n,i-1}_h), w_h \ra, \\
\la a(s^{n,i-1}_h) \bu^{n,i}_h, \bv_h \ra - \la p^{n, i}_h, \nabla \cdot \bv_h \ra + \la \bff_3 (s^{n,i-1}_h), \bv_h \ra &=& 0
\end{eqnarray*}
for all $w_h \in W_h$ and all $\bv_h \in V_h$. Here we use the notation $s^{n,i}_h := s(\Theta^{n,i}_h)$ and, as in the previous section, $s^{n}_h := s(\Theta^{n}_h)$ and $s^{n-1}_h := s(\Theta^{n-1}_h)$. Also, the iteration starts with $\Theta^{n,0}_h := \Theta^{n-1}_h$ and $s^{n,0}_h := s(\Theta^{n-1}_h)$.

We introduce now the errors between two consecutive iterations $i$ and $i-1$:
\begin{displaymath}\begin{array}{rcl@{\,\,\,}rcl}
e_\Theta^{n,i} & = &  \Theta^{n,i}_h - \Theta^{n, i-1}_h, & e_\bq^{n,i}   & = &  \bq^{n,i}_h -  \bq^{n, i-1}_h,\\
e_p^{n,i}     & = & p^{n, i}_h - p^{n, i-1}_h, &e_\bu^{n,i}  & = & \bu^{n,i}_h - \bu^{n, i-1}_h,\\
e_s^{n,i}     & = & s^{n, i}_h - s^{n, i-1}_h := s( \Theta^{n,i}_h )- s(\Theta^{n, i-1}_h).
\end{array}\end{displaymath}
In order to show the convergence of the scheme (\ref{eq1_discrete_lin}) - (\ref{eq4_discrete_lin}) we prove that the iteration is a contraction in $e_\Theta$ and $e_\bq$, together with estimates for $e_p$ and $e_\bu$. The convergence follows by applying the Banach fixed point theorem for $e_\Theta$ and $e_\bq$ and by a similar argument for $e_p$ and $e_\bu$. Note that, in this case, the term involving the factor $L$ will vanish and the limit of the iteration is a solution of Problem $P_h^n$.

The main result, the convergence of the scheme  (\ref{eq1_discrete_lin}) -- (\ref{eq4_discrete_lin}), is stated in the following
\begin{theorem}\label{th_conv_lin} Assuming (A1)-(A5), if the time step $\tau$ is sufficiently small, the linearization scheme  \eqref{eq1_discrete_lin} -\eqref{eq4_discrete_lin} converges linearly.
\end{theorem}
In fact, the time step has to satisfy a mild restriction to guarantee the convergence of the scheme. This restriction is stated in \eqref{restriction_tau} below. In particular, the convergence is robust w.r.t. the mesh size $h$. Moreover, the errors are reduced faster for larger values of $\tau$ than for smaller values, however still satisfying \eqref{restriction_tau}.

The proof of Theorem \ref{th_conv_lin} follows directly from Lemma \ref{lemma1_lin} and Lemma \ref{lemma2_lin} below.

\begin{lemma}\label{lemma1_lin} Let $n \in \N$ be fixed, and $\Theta^{n-1}_h,  p^{n-1}_h \in W_h$ and $\bq^{n-1}_h, \bu^{n-1}_h \in V_h$ be given, solving $P_h^{n-1}$. Further, let
$\Theta^{n,i}_h,  p^{n,i}_h \in W_h$ and $\bq^{n,i}_h, \bu^{n,i}_h \in V_h$ solving $P_h^{n,i}$ for any $i \ge 1, i \in \N$. Assuming (A2)--(A4), there holds
\begin{eqnarray}
\| e_\bu^{n,i} \|^2 & \le & C_1  \| e_s^{n,i-1}\|^2, \label{lemma1_result_1}\\[2ex]
\|  \nabla \cdot e_\bu^{n,i} \|^2 & \le & L^2_{f_2} \| e_s^{n,i-1} \|^2, \label{lemma1_result_2}\\[2ex]
\| e_p^{n,i} \|^2 & \le & C_2 \| e_s^{n,i-1}\|^2, \label{lemma1_result_3}
\end{eqnarray}
where $C_1$ and $C_2$ are two constants not depending on the discretization parameters or the iteration index.
\end{lemma}
\beginproof. The constants $C_1$ and $C_2$ are determined in the proof. With $i > 1$, subtracting 
\eqref{eq3_discrete_lin} and \eqref{eq4_discrete_lin} for $i$ and $i-1$ respectively, one obtains
\begin{eqnarray}
\la \nabla \cdot e_\bu^{n,i}, w_h \ra  &=& \la f_2 (s^{n,i-1}_h) - f_2 (s^{n,i-2}_h), w_h   \ra, \label{lemma1_eq:1}\\
 \la a(s^{n,i-1}_h) \bu^{n,i}_h -  a(s^{n, i-2}_h) \bu^{n, i-1}_h, \bv_h \ra  - \la e_p^{n,i}, \nabla \cdot \bv_h \ra &=&  \la \bff_3(s^{n, i-2}_h) - \bff_3(s^{n,i-1}_h), \bv_h  \ra \label{lemma1_eq:2}
\end{eqnarray}
for all $w_h \in W_h$, $\bv_h \in V_h$. Taking now $w_h = e_p^{n,i} \in W_h$ in \eqref{lemma1_eq:1} and $\bv_h = e_\bu^{n,i} \in V_h$ in \eqref{lemma1_eq:1}, and adding the results one gets
\begin{equation}\label{lemma1_eq:3}
 \la a(s^{n,i-1}_h) \bu^{n,i}_h -  a(s^{n, i-2}_h) \bu^{n, i-1}_h, e_\bu^{n,i}\ra = \la f_2 (s^{n,i-1}_h) - f_2 (s^{n,i-2}_h), e_p^{n,i} \ra + \la \bff_3(s^{n, i-2}_h) - \bff_3(s^{n,i-1}_h),  e_\bu^{n,i}  \ra.
\end{equation}
Using (A2) - (A4), together with Young's inequality, from \eqref{lemma1_eq:3} and for any $\epsilon_1 > 0$ one gets
\begin{equation}\label{lemma1_eq:4}
 \dfrac{a_\star}{2} \| e_\bu^{n,i} \|²^2 \le  \left(\dfrac{M^2_\bu L²^2_a + L^2_{\bff_3}}{a_\star} + \dfrac{L^2_{f_2}}{2 \epsilon_1}\right) \| e^{n,i-1}_s \|^2 + \dfrac{\epsilon_1}{2}  \| e_p^{n,i} \|²^2.
 \end{equation}
Recalling Lemma \ref{lemma_thomas}, a $\bv_h \in V_h$ exists such that $\nabla \cdot \bv_h= e^{n,i}_p $ and $ \| \bv_h \| \le C_{\Omega, d} \|   e^{n,i}_p \|$. Taking this $\bv_h$ as test function in \eqref{lemma1_eq:2}, using (A2)-(A4) gives
\begin{eqnarray}
\|   e^{n,i}_p \|^2 & = & \la a(s^{n,i-1}_h) \bu^{n,i}_h -  a(s^{n, i-2}_h) \bu^{n, i-1}_h, \bv_h \ra + \la \bff_3(s^{n, i-2}_h) - \bff_3(s^{n,i-1}_h), \bv_h  \ra \nonumber \\
& \le &   C_{\Omega, d} \left(  a^\star  \|   e^{n,i}_\bu \|   + (M_\bu L_a + L_{\bff_3}) \|   e^{n,i-1}_s \|  \right)  \|   e^{n,i}_p \|. \nonumber 
\end{eqnarray}
This allows estimating $e^{n,i}_p$ in terms of $e^{n,i}_s$ and $e^{n,i}_\bu$,
\begin{equation}\label{lemma1_eq:5}
\|   e^{n,i}_p \|^2 \le 2 C_{\Omega, d}^2 (a^\star)^2 \|   e^{n,i}_\bu \|^2  + 2 C_{\Omega, d}^2 (M_\bu L_a + L_{\bff_3})^2 \|   e^{n,i-1}_s \|^2.
\end{equation}
With $\epsilon_1 = \dfrac{a_\star}{4  (a^\star)^2 C_{\Omega, d}^2}$, from \eqref{lemma1_eq:4} and \eqref{lemma1_eq:5} one obtains
\begin{equation}\label{lemma1_eq:6}
\dfrac{a_\star}{4} \| e_\bu^{n,i} \|²^2 \le  \left( \dfrac{M^2_\bu L²^2_a + L^2_{\bff_3}}{a_\star} + \dfrac{ 2 (a^\star)^2 C^2_{\Omega, d} L^2_{f_2}}{ a_\star } + \dfrac{a_\star (M_\bu L_a + L_{\bff_3})²^2}{4(a^\star)^2} \right) \| e^{n,i-1}_s \|^2 ,
\end{equation}
which is, in fact \eqref{lemma1_result_1} with $C_1 =  4 \dfrac{M^2_\bu L²^2_a + L^2_{\bff_3} +  2 (a^\star)^2 C^2_{\Omega, d} L^2_{f_2}}{a^2_\star} + \dfrac{(M_\bu L_a + L_{\bff_3})^2}{(a^\star)^2}$. Further, \eqref{lemma1_result_3} follows immediately from \eqref{lemma1_eq:5} and \eqref{lemma1_result_1}, with $C_2 = 2 C^2_{\Omega, d} (a^\star)^2 C_1 + 2C^2_{\Omega, d}(M_\bu L_a + L_{\bff_3})^2$. Finally, \eqref{lemma1_result_2} is a straightforward consequence of \eqref{lemma1_eq:1} and (A3). 
\mbox{  }\finishedproof

\begin{lemma}\label{lemma2_lin} Let $n \in \N$ be fixed, $\Theta^{n-1}_h,  p^{n-1}_h \in W_h$ and $\bq^{n-1}_h, \bu^{n-1}_h \in V_h$ solve ${\bf P_h^{n-1}}$ and  $\Theta^{n,i}_h,  p^{n,i}_h \in W_h$ and $\bq^{n,i}_h, \bu^{n,i}_h \in V_h$ solve ${\bf P_h^{n,i}}$ for any $i \ge 1, i \in \N$.
Assuming (A1)--(A4), it holds
\begin{equation}
\| e_\bq^{n,i} \|^2 \ge  \dfrac{1}{3 C_\Omega^2} \| e_\Theta^{n,i}\|^2 - C_3 \| e_s^{n,i-1}\|^2,\label{lemma2_result_1}
\end{equation}
and
\begin{equation}\label{lemma2_result_2}
\begin{array}{l}
\|  e_\Theta^{n,i} \|^2 + \dfrac{{3C^2_{\Omega, d}} \tau}{2 (3 C^2_{\Omega, d} L +{\tau} )}\|  e_\bq^{n,i} \|^2 \\[2ex] \hspace*{0.5cm}
+ \dfrac{{3 C^2_{\Omega, d}}\left(1 - \tau L (C_3 + 4 L_{\bff_1}^2 + 8 L_{f_w}^2 M^2_\bu + 8 C_1 M_{f_w}^2) \right )}{L({3 C^2_{\Omega, d}} L +{\tau})} \|  e_s^{n,i-1} \|^2
\le  \dfrac{{3  C^2_{\Omega, d}} L }{{3 C^2_{\Omega, d}}L  +{\tau} } \| e_\Theta^{n,i-i} \|^2,
\end{array}
\end{equation}
where $C_3 = M^2_{f_w} C_1 + (M_{f_w} L_{f_w}  + L_{\bff_1})^2$ is not depending on the discretization parameters.
\end{lemma}
\beginproof. 
Subtracting \eqref{eq1_discrete_lin} and \eqref{eq2_discrete_lin} for $i$ and $i-1$ respectively gives
\begin{eqnarray}
\la L (e_\Theta^{n,i} - e_\Theta^{n,i-1}) + e_s^{n,i-1}, w_h \ra  + \tau \la \nabla \cdot e_\bq^{n,i}, w_h \ra &=& 0,\label{lemma2_eq:1}\\[1ex]
\la e_\bq^{n,i}, \bv_h \ra -\la e_\Theta^{n,i},  \nabla  \cdot \bv_h \ra -  \la f_w (s^{n,i-1}_h) \bu^{n,i}_h -f_w (s^{n,i-2}_h) \bu^{n, i-1}_h ,  \bv_h \ra &=& \la \bff_1(s^{n,i-1}_h) - \bff_1(s^{n, i-2}_h), \bv_h \ra. \nonumber \\\label{lemma2_eq:2}
\end{eqnarray}
By Lemma \ref{lemma_thomas}, there exists a $\bv_h \in V_h $ such that $\nabla \cdot\bv_h =  e_\Theta^{n,i}$ and $\| \bv_h \| \le C_{\Omega, d} \| e_\Theta^{n,i} \| $. Taking this  $\bv_h $ as test function in \eqref{lemma2_eq:2} and using (A3) and (A4) gives
\begin{eqnarray}
 \| e_\Theta^{n,i} \|^2 & = &  \la e_\bq^{n,i}, \bv_h \ra + \la f_w (s^{n,i-1}_h) \bu^{n,i}_h - f_w (s^{n, i-2}_h) \bu^{n, i-1}_h ,  \bv_h \ra + \la \bff_1(s^{n,i-1}_h) - \bff_1(s^{n, i-2}_h) , \bv_h \ra\nonumber \\
& \le &  \left( \| e_\bq^{n,i} \|  + \| f_w (s^{n,i-1}_h) e_\bu^{n,i}\| + \| (f_w (s^{n,i-1}_h) -  f_w (s^{n, i-2}_h) ) \bu^{n, i-1}_h\|   + L_{\bff_1} \| e_s^{n,i-1} \| \right) \| \bv_h \|  \nonumber \\
& \le & \left(  \| e_\bq^{n,i} \|  + M_{f_w} \| e_\bu^{n,i} \|   + (M_{\bu} L_{f_w}  + L_{\bff_1}) \| e_s^{n,i-1} \|   \right) C_{\Omega, d} \| e_\Theta^{n,i} \|.  \label{lemma2_eq:3}
\end{eqnarray}
Combining \eqref{lemma2_eq:3} with \eqref{lemma1_result_1} further implies
\begin{equation}\label{lemma2_eq:4}
 \| e_\Theta^{n,i} \|^2 \le 3 C^2_{\Omega, d}²  \| e_\bq^{n,i} \|^2 + 3 C^2_{\Omega, d} C_3 \| e_s^{n,i-1} \|^2.
\end{equation}
where $C_3 =  M^2_{f_w} C_1 + (M_{\bu} L_{f_w}  + L_{\bff_1})^2.$ From \eqref{lemma2_eq:4}, \eqref{lemma2_result_1} follows immediately.

To prove \eqref{lemma2_result_2}, one takes $w_h = e_\Theta^{n,i} \in W_h$ in \eqref{lemma2_eq:1} and $\bv_h = \tau e_\bq^{n,i}$ in \eqref{lemma2_eq:2}, add the resulting and obtains
\begin{displaymath}\label{lemma2_eq:5}
\begin{array}{l}
\ds L \la (e_\Theta^{n,i} - e_\Theta^{n,i-1}) + e_s^{n,i-1}, e_\Theta^{n,i} \ra  + \tau \| e_\bq^{n,i} \|²^2 = \\[1ex]
\ds \hspace*{3cm} \tau \la f_w (s^{n,i-1}_h) \bu^{n,i}_h -f_w (s^{n, i-2}_h) \bu^{n, i-2}_h,  e_\bq^{n,i} \ra + \tau \la \bff_1(s^{n,i-1}_h) - \bff_1(s^{n, i-2}_h) , e_\bq^{n,i}\ra.
\end{array}
\end{displaymath}
This further implies
\begin{equation}\label{lemma2_eq:6}
\begin{array}{l}
\ds \dfrac{L}{2} \| e_\Theta^{n,i}\|^2 + \dfrac{L}{2} \| e_\Theta^{n,i}- e_\Theta^{n,i-1}\|^2 + \la e_s^{n,i-1}, e_\Theta^{n,i-1} \ra  + \tau \| e_\bq^{n,i} \|²^2 = \dfrac{L}{2} \| e_\Theta^{n,i-1}\|^2 + \la e_s^{n,i-1}, e_\Theta^{n,i-1} -  e_\Theta^{n,i}\ra  \\[2ex]
\ds  \hspace*{0.5cm} + \tau \la f_w (s^{n,i-1}_h) \bu^{n,i}_h -f_w (s^{n, i-2}_h) \bu^{n, i-2}_h,  e_\bq^{n,i} \ra + \tau \la \bff_1(s^{n,i-1}_h) - \bff_1(s^{n, i-2}_h) , e_\bq^{n,i} \ra.
\end{array}
\end{equation}
By the monotonicity and the Lipschitz continuity of $s(\cdot)$ as stated in (A1) there holds
\begin{equation}\label{lemma2_eq:7}
\la e_s^{n,i-1}, e_\Theta^{n,i-1} \ra \ge \dfrac{1}{L_s} \| e_s^{n,i-1}\|^2 \ge \dfrac{1}{L} \| e_s^{n,i-1}\|^2.
\end{equation}
From \eqref{lemma2_result_1} and \eqref{lemma2_eq:7}, by (A1) - (A4) and Young's inequality, \eqref{lemma2_eq:6} implies
\begin{equation}\label{lemma2_eq:8}
\begin{array}{l}
\ds (\dfrac{L}{2} + \dfrac{\tau}{6 C^2_{\Omega,d}}) \| e_\Theta^{n,i}\|^2 + \dfrac{L}{2} \| e_\Theta^{n,i}- e_\Theta^{n,i-1}\|^2 +  (\dfrac{1}{L} - \tau \dfrac{C_3}{2})\| e_s^{n,i-1}\|^2  + \dfrac{\tau}{2} \| e_\bq^{n,i} \|²^2 \\[1ex]
\hspace*{1cm} \le  \dfrac{L}{2} \| e_\Theta^{n,i-1}\|^2  +  \dfrac{1}{2 L }\| e_s^{n,i-1}\|^2 + \dfrac{L}{2} \| e_\Theta^{n,i}- e_\Theta^{n,i-1}\|^2 + 2 \tau L_{\bff_1}^2 \| e_s^{n,i-1}\|^2   \\[1ex]
\hspace*{1cm} + \dfrac{\tau}{8} \| e_\bq^{n,i}\|^2  + 4 \tau L_{f_w}^2 M^2_\bu  \| e_s^{n,i-1}\|^2 + 4 \tau M_{f_w}^2 \| e_\bu^{n,i} \|^2 + \dfrac{\tau}{8} \| e_\bq^{n,i} \|^2 .
\end{array}
\end{equation}
This rewrites as
\begin{equation}\label{lemma2_eq:9}
\begin{array}{l}
\ds (\dfrac{L}{2} + \dfrac{\tau}{6 C^2_{\Omega,d}}) \| e_\Theta^{n,i}\|^2 +  (\dfrac{1}{L} - \tau \dfrac{C_3}{2})\| e_s^{n,i-1}\|^2  + \dfrac{\tau}{4} \| e_\bq^{n,i} \|^2 \\[1ex]
\hspace*{1cm} \le  \dfrac{L}{2} \| e_\Theta^{n,i-1}\|^2  +  (\dfrac{1}{2 L } + 2 \tau L_{\bff_1}^2 + 4 \tau L_{f_w}^2 M^2_\bu) \| e_s^{n,i-1}\|^2 + 4 \tau M_{f_w}^2 \| e_\bu^{n,i} \|^2.
\end{array}
\end{equation}
Using now \eqref{lemma1_result_1} in \eqref{lemma2_eq:9} and rearranging the terms leads to
\begin{equation}\label{lemma2_eq:10}
\begin{array}{l}
\ds (\dfrac{L}{2} + \dfrac{\tau}{6 C^2_{\Omega,d}}) \| e_\Theta^{n,i}\|^2 +  \left(\dfrac{1}{2L} - \tau (\dfrac{C_3}{2} +  2 L_{\bff_1}^2 + 4 L_{f_w}^2 M^2_\bu + 4 C_1 M_{f_w}^2) \right )\| e_s^{n,i-1}\|^2 \\
 \hspace*{1cm}+ \dfrac{\tau}{4} \| e_\bq^{n,i} \|²^2  \le  \dfrac{L}{2} \| e_\Theta^{n,i-1}\|^2,
\end{array}
\end{equation}
which is nothing else as the result \eqref{lemma2_result_2}.
\finishedproof

\begin{remark} The estimate \eqref{lemma2_result_2} is not practical unless the factor multiplying the last term on the left is positive. This gives a restriction on the time step,
\begin{equation}\label{restriction_tau}
\tau \leq \dfrac{1}{L \left(C_3 + 4 L_{\bff_1}^2 + 8 L_{f_w}^2 M^2_\bu + 8 C_1 M_{f_w}^2\right)} ,
\end{equation}
where $C_1 =  4 \frac{M^2_\bu L²^2_a + L^2_{\bff_3} +  (a^\star)^2 C^2_{\Omega, d} L^2_{f_2}}{a^2_\star} + \frac{2 (M_\bu L_a + L_{\bff_3})^2}{ (a^\star)^2} $ and $C_3 =  M^2_{f_w} C_1 + (M_{\bu} L_{f_w}  + L_{\bff_1})^2$. This is a mild condition because it does not depend on the grid size. In this sense, it is superior to the conditions guaranteeing the stability of an explicit scheme in time, or to the typical conditions guaranteeing the convergence of the Newton method for degenerate parabolic problems (see e.g. \cite{radu2006, radu2011}).
\end{remark}

\begin{remark} The constant $L$ is independent on the discretization parameters, but must be chosen greater than the Lipschitz constant of the function $s(\cdot)$, i.e. $L_s$.
\end{remark}

\begin{remark} One can use the $L$-scheme  \eqref{eq1_discrete_lin}--\eqref{eq4_discrete_lin} also in combination with the Newton method. The goal is to combine the robustness of the $L$-scheme, which converges regardless of the starting point, with the quadratic convergence of the Newton Method, which requires instead a starting point close to the solution. Specifically, at each time step one can perform a few $L$-scheme iterations, followed by Newton iterations. In this way one enhances the robustness of the Newton method and reduces the severe restriction on the time step that guarantees its convergence. We refer to \cite{list}, where this strategy is studied for solving the Richards equation. Still for the Richards equation, A similar idea was also proposed in \cite{ackerer}, but there the modified Picard method was used to improve the robustness of the Newton method.
\end{remark}

\begin{remark} We refer to \cite{radu2016hoelder} for a $L$-scheme applied to the case of a  H\"older continuous saturation.
\end{remark}

\section{Numerical results}\label{sec:simulation}
In this section we present two numerical studies, one concentrating on the convergence of the backward Euler/MFEM discretization and one on the convergence of the linearization scheme. For more numerical examples we refer to \cite{radu_NumMath, radu_SIAM} (for convergence of the discretization error) and \cite{list} (for linearization schemes). These papers are considering Richards' equation (which is just a particular case of the two-phase flow considered in this paper). We further refer to \cite{radu2015} for an example concerning the linearization method for two-phase flow and MPFA.

We consider here two problems. The first is defined in a two-dimensional domain $\Omega = (0,1) \times (0,1)$ and has the analytical solution given in \eqref{eq:manufactured}. For this, a source term was added to \eqref{eq_classic_1}:
\begin{align*}
f (t,x,y) = 2tx^2 (1-x)^2 y^2 (1-y)^2 + 2t x(1-x) + 2 t y (1-y)\\
+ t^2 y^3 (1-y)^3 \left ( 10 x^4 - 20 x^3 +12 x^2 - 2x \right ) \\
+ t^2 x^3 (1-x)^3 \left ( 10 y^4 - 20 y^3 +12 y^2 - 2y \right ),
\end{align*}
and we choose appropriate initial and (Dirichlet) boundary conditions. For the spatial discretization, we use a rectangular and uniform mesh, whereas for the time discretization we choose a uniform time step with final time $T = 1$. In accordance with the estimate in Theorem \ref{th_main_result_estimates}, we will consider a sequence of discretizations with halving the spatial mesh size $h$ and reducing the time step $\tau$ one-fourth.

Recalling the system of equations \eqref{eq_classic_1} - \eqref{eq_classic_4}, the solution and coefficient functions are given by
\begin{align} \label{eq:manufactured}
\begin{array} {lllll}
p =  x (1-x) y(1-y), & \; \Theta = t x (1-x) y (1-y), & \; s(\Theta)= \Theta^2, \;  \lambda_w = s, \\
 \lambda_o = 1-s, \; & \textbf{f}_1 = 0, \; & f_2 = 2x (1-x) + 2 y(1-y), \; & \textbf{f}_3 = 0.\\
\end{array}
\end{align}

\begin{table}[!h]
{\footnotesize
\begin{center}
  \begin{tabular}{ | c | c | c | c | c | c | c | c | c | c |}
    \hline
    $h$ & $\tau$ & $E_p$  &   conv rate           &        $E_{s\Theta}$  &       conv rate  &  $E_{\Theta}$  & conv rate      & $E_s$  &   conv rate        \\ [1em] \hline
    $\frac{1}{4}$ & $\frac{1}{5}$   & $1.96 E - 4 $ &   {-- }  &         3.35E-6     &    {--}       &     9.14 E-5          & {--}      & $1.82 E - 7$ & {--}     \\ [1em] \hline
    $\frac{1}{8}$ & $\frac{1}{20}$  & $4.48 E-5$ & 2.13   &         5.53E-7      &   2.59         &   1.66E-5         & 2.46    & $2.88 E-8$ & 2.66       \\ [1em] \hline
    $\frac{1}{16}$ & $\frac{1}{80}$  & $1.11 E-5$ & 2.01   &        1.28 E-7   &    2.11          &  3.9 E-6         & 2.09       & $6.67 E-9$ & 2.11      \\ [1em] \hline
    $ \frac{1}{32}$ & $\frac{1}{320}$  & $2.72 E-6$ & 2.00 &     3.12E -8 &   2.04          & 9.53 E-7       & 2.03            & $1.61 E - 9$ & 2.05     \\ [1em] \hline
  \end{tabular}
  \caption{Convergence rates for the manufactured solution \eqref{eq:manufactured}.}   \label{fig:convtable}
\end{center}
}\end{table}

The results are presented in Table \ref{fig:convtable}, with the errors given by:
\begin{align*}
\begin{array}{lllll}
E_{p} =  \sum_{n=1}^N \tau \| \opn -  p^n_h \| ^2 , & E_{\Theta} = \sum_{n=1}^N \int_{t_{n-1}}^{t_{n}} \|  \Theta(t) - \Theta^n_h \|^2 \, dt,  \\[1em]
E_{s} = \sum_{n=1}^N \int_{t_{n-1}}^{t_{n}} \|  s(t) - s^n_h \|^2 \, dt, &   E_{s\Theta} = \sum_{n=1}^N \int_{t_{n-1}}^{t_{n}} \la s(\Theta(t)) - s(\Theta^n_h), \Theta(t) - \Theta^n_h \ra \, dt. \\[1em]
\end{array}
\end{align*}
The convergence rate is computed by
\[ \text{conv\;rate} (i) = \dfrac{ \log E_{z} (i+1) - \log E_{z} (i) }{\log h(i+1) - \log h(i)},\]
where $i$ is the array index and $z$ stands for either $p$, $\Theta$, $s$, or $s\Theta$, as shown in Table \ref{fig:convtable}. We see that the convergence rate is $2$ which is as expected.

Next, we discuss another example where we consider 3D rectangular grids of different sizes and study the convergence of linearization scheme. The computational domain is now the unit cube. We use the following constitutive relationships
\begin{align*}
k_{rw} = s^2, \;  k_{ro} = (1-s)^2, \; \Theta = \sqrt{s},
\end{align*}
and use the following parameters
\begin{align*}
T = 50 \text{ days }, \; \tau = 0.5 \text{ day },  \; L = 2,   k = 10^{-6} \; \mathrm{ m}^2,  \; \mu_w = 1 \text{ cP} \;, \mu_o = 10 \text{ cP}.
\end{align*}
In terms of the coefficients in the model equations \eqref{eq_classic_1} - \eqref{eq_classic_4}, these choices correspond to:
\begin{align*}
a(s) = 10^{-6}
 \frac{1}{s^2 + (1-s)^2}, \; f_w(s) = \frac{s^2}{s^2 + (1-s)^2}, \;  \bf{f}_1 = 0,\; \bf{f}_3 = 0.
\end{align*}

\begin{figure}[h!]
   \includegraphics[width=\textwidth]{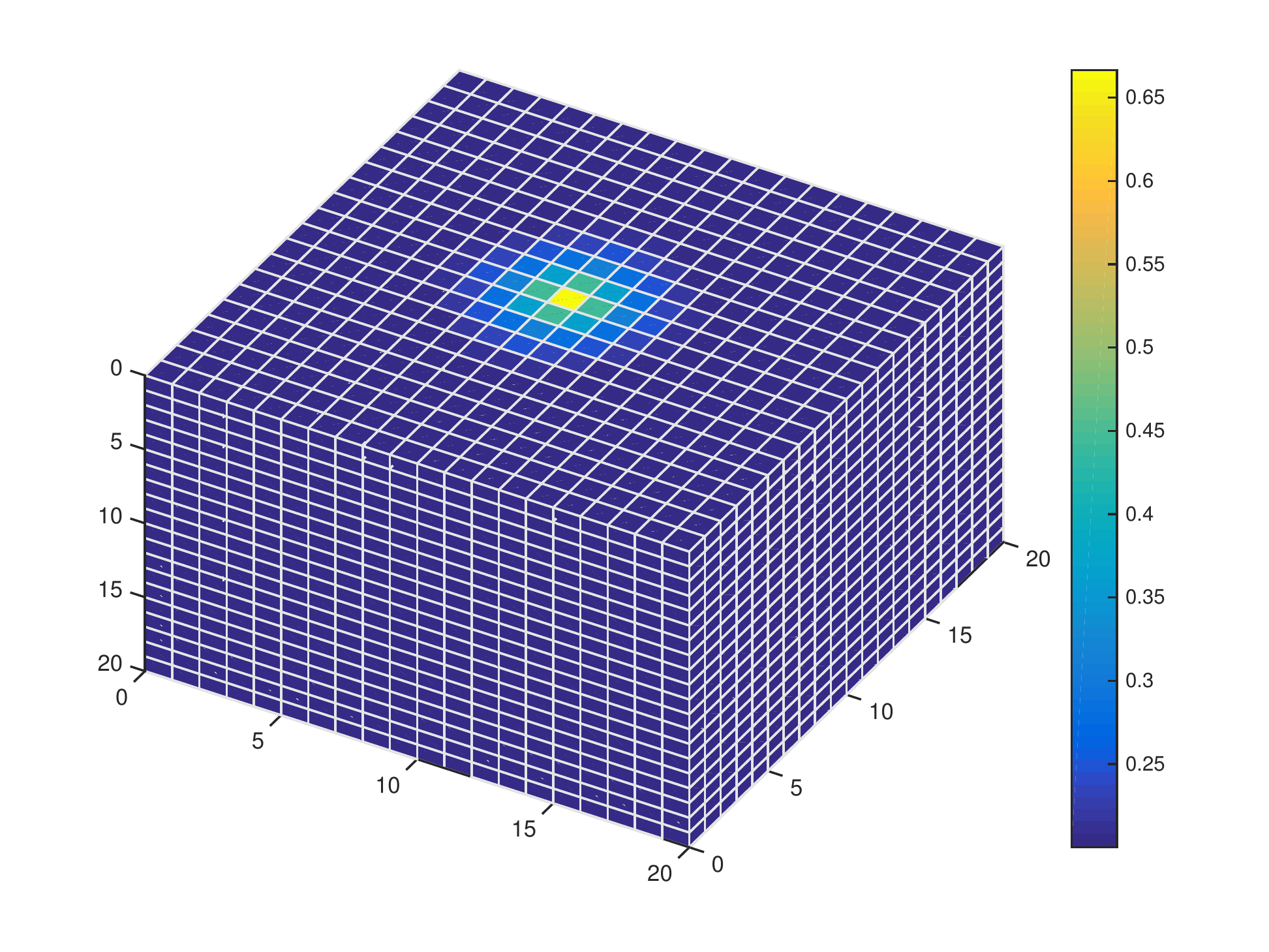}
  \caption{Plot of saturation at time $T = 20$ days. In the middle of the grid, water is injected at a constant rate.} \label{fig:saturationplot}
\end{figure}

\begin{figure}[h!]
      \includegraphics[width=\textwidth]{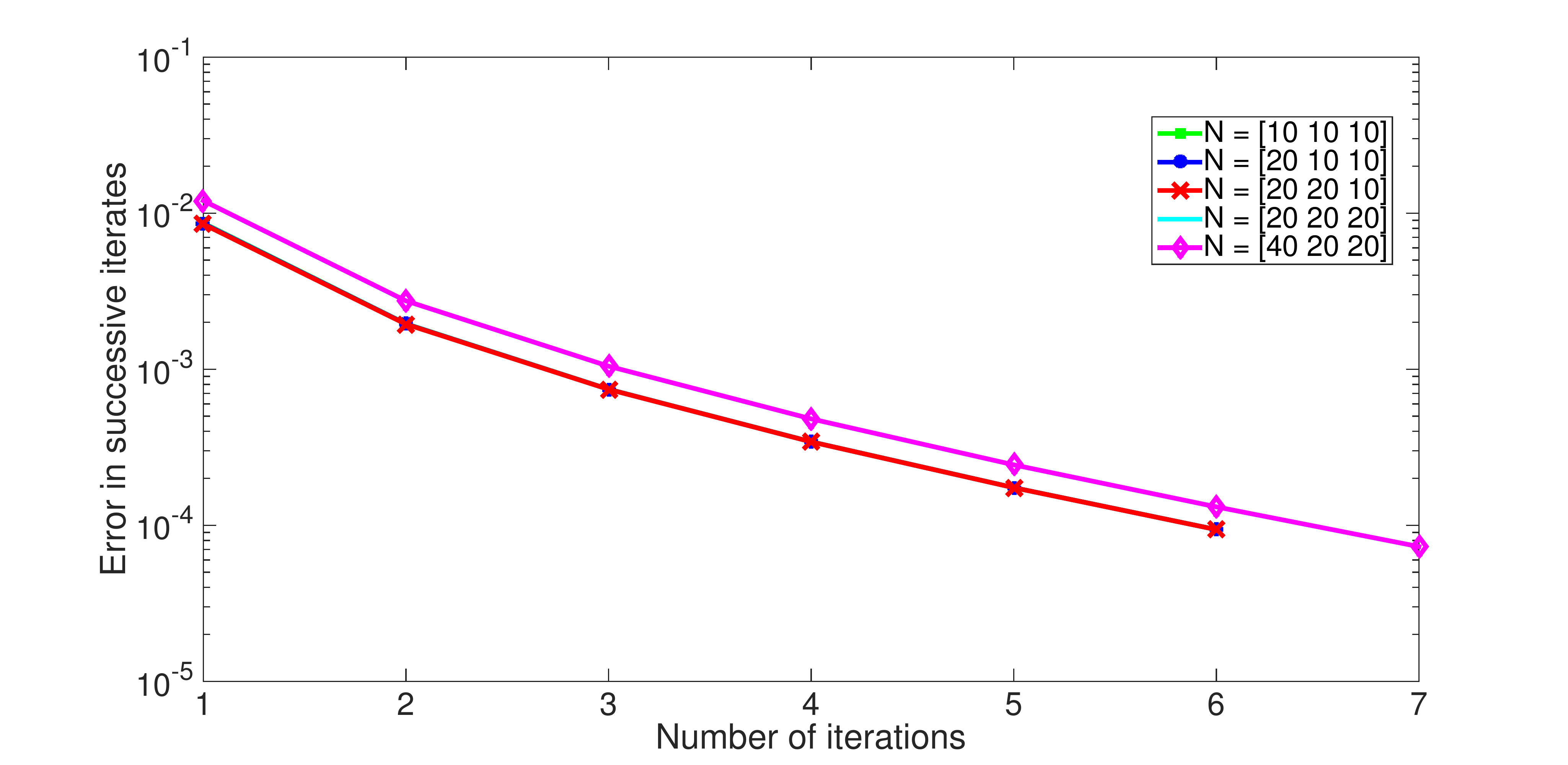}
  \caption{Convergence of the $L$-scheme for different grids and a Lipschitz continuous $s(\cdot)$. The legend in the figure shows the grid sizes which have been taken. } \label{fig:iterationvsh}
\end{figure}

\begin{figure}[h!]
      \includegraphics[width=\textwidth]{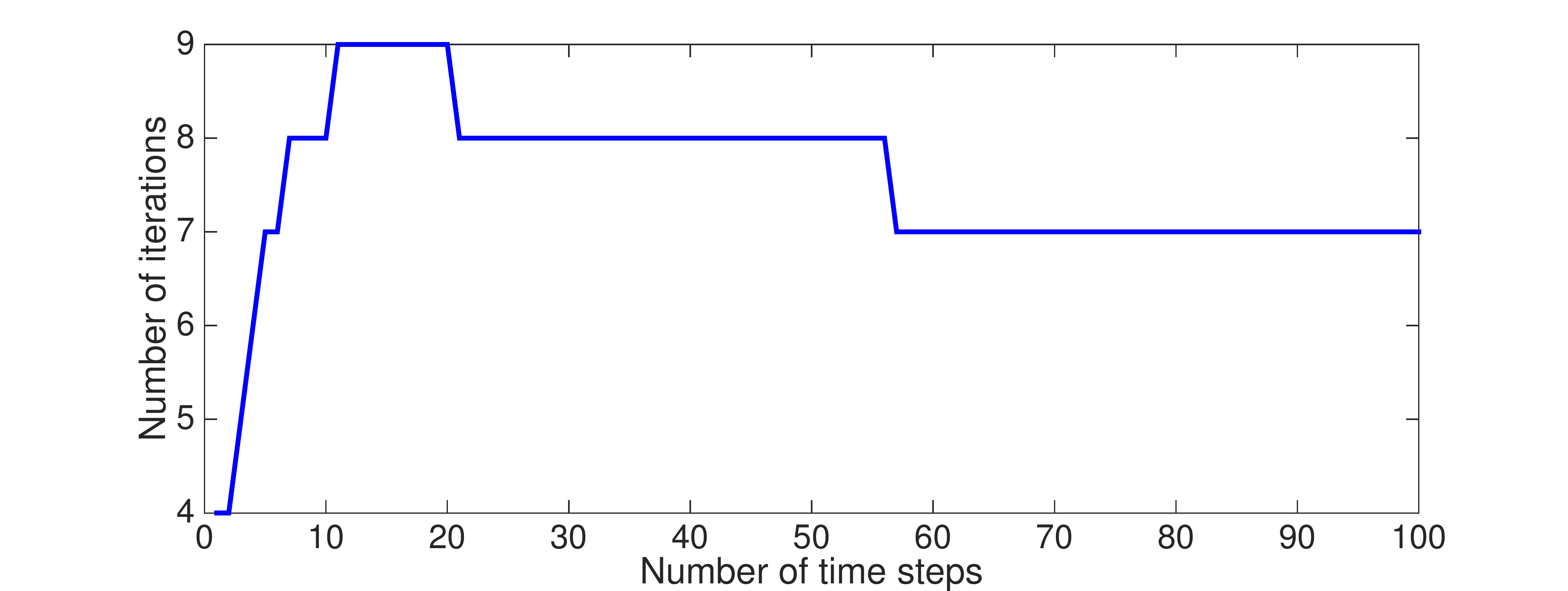}
   \caption{Number of iterations versus time step.}   \label{fig:iterationvst}
\end{figure}

For the pressure, we use Dirichlet boundary conditions at the left ($p =0 \text{ at } x = 0$) and right sides ($p = 10 \text{ at } x = 1$) and homogeneous Neumann at the rest of the boundaries.   For the saturation, we use no flow boundary conditions and consider an injection at the center of the cells $f_2 = 10^{-5}$ m$^3$/s. For the grid of size $nx = 20, ny=  20, nz= 20$, the saturation plot at $T = 20$ days is shown in Figure \ref{fig:saturationplot}. In Figure \ref{fig:iterationvsh}, we show the convergence of linear iteration in one time step (at $T = 20$ days) for different grid sizes. We see that the convergence is rather independent of the problem size. Moreover, as shown in Figure \ref{fig:iterationvst} we show that the number of linear iterations is not very sensitive (5-9 iterations) for the given problem at any time step. 

\section{Conclusions}\label{sec:conclusions} We considered a mathematical model for two-phase flow in porous media. The model was formulated in terms of a global and a complementary pressure. A fully implicit, mass conservative numerical scheme was proposed for solving it numerically. The scheme is based on backward Euler for the discretization in time and mixed finite elements (lowest order Raviart - Thomas elements) for the spatial discretization. The scheme was shown to be convergent. Moreover, order of convergence estimates are provided. For solving the non-linear systems at each time step we considered a robust, first order convergent linearization method, called the  $L$-method. The convergence of the linearization scheme was rigorously shown. The $L$-method does not involve the computation of any derivatives, and does not require any regularization step. The method can be used also to enhance the robustness of Newton's method, as was done in \cite{list}. The convergence rate of the method does not depend on the mesh diameter. Numerical examples have been shown to sustain the theoretical results.\\[2ex]

{\bf Acknowledgments.}
The authors are members of the International Research Training Group  NUPUS funded by the German Research Foundation DFG (GRK 1398), the  Netherlands Organisation for Scientific Research NWO (DN 81-754) and  by the Research Council of Norway (215627). F. A. Radu and I. S. Pop acknowledge the support of Statoil through the Akademia agreement. J.M. Nordbotten acknoweldges the NWO support through the Visitors Grant 040.11.351. Part of this work was done during F.A. Radu's sabbatical in Eindhoven, we acknowledge for this the support of Meltzer foundation, of University of Bergen and the NWO Visitors Grant 040.11.499.


\end{document}